\documentclass[12pt]{article}
\usepackage[active]{srcltx}
\usepackage{xcolor}
\usepackage{amscd,amsfonts,amssymb,amscd,amsmath,latexsym}
\usepackage[pdftex]{hyperref}
\usepackage{setspace}
\usepackage[normalem]{ulem}
\usepackage[all]{xy}
\usepackage{mathrsfs}
\usepackage{graphicx}
\usepackage{amsthm}

\title{T-Duality via Gerby Geometry and  Reductions}
 \author{Ulrich Bunke and Thomas Nikolaus\thanks{{Fakult\"at f\"ur Mathematik,
Universit{\"a}t Regensburg,
93040 Regensburg,
{Germany}, \href{mailto:ulrich.bunke@mathematik.uni-regensburg.de}{ulrich.bunke@mathematik.uni-regensburg.de}\ ,
\href{mailto:thomas1.nikolaus@mathematik.uni-regensburg.de}{thomas1.nikolaus@mathematik.uni-regensburg.de}\ .
}}
}

\newtheorem{theorem}{Theorem}[section] 
\newtheorem{prop}[theorem]{Proposition}
\newtheorem{lem}[theorem]{Lemma}
\newtheorem{ddd}[theorem]{Definition}
\newtheorem{kor}[theorem]{Corollary}

\newtheorem{fact}[theorem]{Fact}

\theoremstyle{remark}

\theoremstyle{definition}
\newtheorem{rem2}[theorem]{Remark}

\newcommand{\PSh}{\mathrm{PSh}}
\newcommand{\taut}{  univ }

\newcommand{\colim}{{\tt colim}}

\newcommand{\Z}{\mathbb{Z}}

\renewcommand{\proof}{{\it Proof.$\:\:\:\:$}}

\newcommand{\kaaa}{{\frak k}}

\newcommand{\taaa}{{\frak t}}

\newcommand{\R}{\mathbb{R}}

\newcommand{\gaaa}{{\frak g}}

\newcommand{\Aut}{{\tt Aut}}

\newcommand{\cV}{\mathcal{V}}

\newcommand{\cI}{\mathcal{I}}

\newcommand{\cU}{\mathcal{U}}
\newcommand{\Hom}{{\tt Hom}}

\newcommand{\Ext}{{\tt Ext}}

\newcommand{\im}{{\tt im}}
\newcommand{\sign}{{\tt sign}}

\newcommand{\cF}{\mathcal{F}}

\newcommand{\Ad}{{\tt  Ad}}

\newcommand{\Gr}{{\tt Gr}}

\newcommand{\id}{{\tt id}}

\newcommand{\nat}{\mathbb{N}}
\newcommand{\supp}{{\tt supp}}

\newcommand{\cB}{\mathcal{B}}
\def\imath{{i}}

\def\hB{\hspace*{\fill}$\Box$ \newline\noindent}

\def\hB{\hspace*{\fill}$\Box$ \\[0.5cm]\noindent}
\newcommand{\cL}{\mathcal{L}}
 \newcommand{\cG}{\mathcal{G}}

\newcommand{\pr}{{\tt pr}}

\newcommand{\Ab}{{\mathrm{Ab}}}

\newcommand{\Mf}{{\mathbf{ Mf}}}

\newcommand{\bbG}{\underline{S}^{1}}

\newcommand{\tho}[1]{\textcolor{blue}{#1}}

\begin{document}
\maketitle
\begin{abstract}
We consider  {topological} $T$-duality of torus bundles equipped with $\bbG$-gerbes. 
We show how a geometry on the gerbe determines a reduction of its band to the subsheaf of $S^{1}$-valued functions which 
are constant along the torus fibres. We observe that such a reduction is exactly the additional datum 
needed for {the} construction of {a} $T$-dual pair. We illustrate the theory by working out the example of the canonical lifting gerbe on a compact Lie group
{which is a torus bundles over the associated flag manifold. It was a recent observation of Daenzer and van Erp \cite{2012arXiv1211.0763D} that for certain compact Lie groups and a particular choice of the gerbe, the $T$-dual torus bundle is 
given by the Langlands dual group.} 
 \end{abstract}

\tableofcontents
 
\section{Introduction}
 
Topological $T$-duality is a relation between torus bundles equipped with a twist. 
In the present paper we determine which choice of geometric data can be used to define $T$-duals in a functorial manner.

We work in stacks on smooth manifolds and realize twists as gerbes with band $\bbG$, the sheaf of smooth $S^{1}$-valued functions. 
Over the total space of a torus bundle the sheaf $\bbG$ has the subsheaf of functions which are constant along the fibres.
We observe that a reduction of the band of the twist to this subsheaf
is exactly the datum needed to define a $T$-dual. Furthermore, the gerbe with the reduced band must be trivializable along the fibres.
In the present paper we {analyse} when a geometry on the original gerbe induces such a reduction, and under which additional conditions 
this reduction is fibrewise trivializable. As an example we work out in detail the case of lifting gerbes and consider in particular 
the case of compact Lie groups as total spaces of torus bundles over their flag manifolds. 

\bigskip

{In order to put the paper into a context, we start this introduction with a review of $T$-duality, and then explain our results in greater detail.}
Consider an $n$-dimensional torus $T$,  a $T$-principal bundle $E\to B$ and a twist $H\to E$.
Depending on the set-up, twists are {realized} alternatively by
\begin{enumerate}
\item  a bundle $H\to E$ with fibre $K(\Z,2)$ classified by a map $E\to K(\Z,3)$, where $K(\Z,n)$ is an Eilenberg-MacLane space,
\item gerbes $H\to E$ with band $\bbG$, the sheaf of $S^{1}$-valued functions,
\item bundles $H\to E$ of algebras of compact operators.
\end{enumerate}
In all cases the equivalence class of a twist $H\to E$  is determined   by a  class $[H\to E]\in H^{3}(E;\Z)$.
The datum $$H\to E\to B$$ is called a pair. 

Topological $T$-duality is a relation between two such pairs $$H\to E\to B\ ,\quad \hat H\to \hat E\to B$$ over the same base $B$. We use the term {``}topological'' since 
only the topology of these objects matters. Topological $T$-duality models the underlying topology of a geometric $T$-duality. The term $T$-duality is motivated from string theory where a pair with additional geometric structures is considered as a background for a certain type {of} string theory. Its $T$-dual considered as a background of a different string theory should produce a physically equivalent theory.
These string theoretic considerations predict isomorphisms between topological invariants of the $T$-dual pairs $H\to E\to B$ and $\hat H\to \hat E\to B$, most prominently between their twisted real cohomology and twisted $K$-theory. We refer to the literature cited below for more information on the physical background.

\bigskip

The first instance of topological $T$-duality has been considered {in}
\cite{MR2080959}, \cite{MR2116165},  \cite{MR2062361}, \cite{MR2202291}.
In these papers the concept of  $T$-duality was based on a differential form level  consideration of the 
$T$-duality isomorphism for the twisted de Rham cohomology. In particular it  has been observed for the first time that the presence of the twist $H\to E$ influences the  topology of the dual $T$-bundle $\hat E\to B$.

In order to establish the $T$-duality isomorphism in twisted $K$-theory in
\cite{MR2130624}, \cite{MR2287642}, \cite{MR2797285} {the first author and coauthors} introduced, in increasing generality, a concept of $T$-duality based on the notion of a $T$-duality triple.
A $T$-duality triple is a diagram of the form
$$\xymatrix{ &H\times_{B}\hat E\ar[dl]\ar[dr]&&E\times_{B}\hat H\ar[dr]\ar[dl]\ar[ll]^{u}&\\H\ar[dr]&&E\times_{B}\hat E\ar[dr]\ar[dl]&&\hat H\ar[dl]\\&E\ar[dr]&&\hat E\ar[dl]&\\&&B&&}$$
which includes the two pairs on the left- and right-hand sides, and where the map $u$ of $\bbG$-banded gerbes must satisfy a {`Poincar\'e bundle condition'}.
For details we refer to  \cite[Definition 2.8]{MR2287642},
\cite[Definition 4.1.3]{MR2797285}.
\begin{ddd}
A pair $\hat H\to \hat E\to B$ is topologically $T$-dual to a pair $H\to E\to B$ if  it fits into a $T$-duality diagram as above.
\end{ddd}
Given the pair $H\to E\to H$ the natural {problem} is now to classify its extensions to $T$-duality triples and the possible dual pairs. If $n=1$, then there is always a uniquely determined dual, but for $n\ge 2$
there are obstructions against the extension to a $T$-duality triple, and if such extensions exists, they depend on additional choices. The classification problem is solved in \cite{MR2287642}.

There is an alternative approach to topological $T$-duality via non-commutative topology 
\cite{MR2116734},\cite{MR2222224}, \cite{MR2327179}. In this case one {realizes} the twist $H\to E$ as {a} bundle of $C^{*}$-algebras and considers $E$ as the prime ideal spectrum of the algebra of sections $C(E,H)$.
 Let $\tilde T\to T$ be the universal covering of the torus. In order to define a $T$-dual one must lift the $T$-action on $E$ to a $\tilde T$-action on $C(E,H)$.  The $T$-dual, if it exists, is then defined as the prime ideal spectrum of the cross product, i.e. by  $C(\hat E,\hat H)\cong C(E,H)\rtimes {\tilde T}$.

 It has been verified in the thesis of A. Schneider \cite{2007arXiv0712.0260S} that the classification 
 of classical $T$-duals via non-commutative topology is equivalent to the classification via $T$-duality triples. Advantages of the non-commutative topology approach are that it yields the $T$-duality isomorphism in twisted $K$-theory naturally as Connes' Thom isomorphism and also
produces exotic $T$-duals in non-commutative topology if classical $T$-duals do not exist.
Moreover, the cross product gives an explicit {and functorial} construction of $T$-duals.

\bigskip

One can consider {extensions} of the theory which  incorporate additional  automorphisms, group actions,  
or torus actions with fixed points, {see} \cite{MR2246781}, \cite{MR2369414}, \cite{MR2491618}, \cite{2009arXiv0906.3734S}, \cite{2012arXiv1211.2890P}.  
These papers mostly consider the non-commutative topology approach. In the topological context much of this theory can be subsumed by working in the context of topological stacks  as in \cite{MR2797285}.

\bigskip

Let us now consider the question which additional structures on the pair $H\to E\to B$ 
determine a $T$-dual or its extension to a $T$-duality triple. In \cite{MR2287642} (see also Theorem \ref{feb2001}) we explain the choice of cohomological data which characterize the $T$-dual up to equivalence.  But this data is not sufficient for an actual construction of the $T$-dual.

\bigskip

In the non-commutative topology approach the choice of a lift of the $T$-action on $E$ to a $\tilde T$-action on $C(E,H)$ allows an actual construction of a $T$-dual.  Similarly, in the topological context   the choice of the structure of a torsor on $H\to B$ over a  Picard stack  
allows to construct the $T$-dual by a version of {Pontrjagin} duality for locally compact Picard stacks  \cite{MR2482327}. This was motivated by    analogous ideas  in \cite{MR2399730} in the context of algebraic geometry. 
 
\bigskip

In the present paper we describe a very natural choice of data which allows for a canonical construction of the $T$-dual. Recall that $H\to E$ is a gerbe with band $\bbG$.
This sheaf contains the subsheaf $\bbG_{sm/B}$ of $S^{1}$-valued functions which are locally constant along the fibres of $E\to B$. 
The   {additional} data is the choice of a {$\bbG_{sm/B}$}-reduction $H^{\flat}\to E$ of the 
gerbe $H\to E$ (Definition \ref{feb1201}). The $T$-dual $\hat H$ then simply appears as the stack of fibrewise 
trivializations of $H^{\flat}\to E$, and $\hat E$ is the moduli space of $\hat H$ (Definition \ref{feb1903}).
In particular, a $T$-dual exists if and only if the gerbe $H\to E$ admits a reduction which is fibrewise trivializable.
{A more detailed explanation will be given in Section \ref{feb2002}}.
This constructive definition of the $T$-dual is functorial in the data. {Therefore it allows to transport additional structures like automorphisms or 
geometric data from $H\to E\to B$ to the dual side $\hat H\to \hat E\to B$.}
The idea that the $T$-dual is related to the fibrewise trivializations of the original data
has been around for a while but we did not know how to make it precise in the differential-topological context.

\bigskip

In Section \ref{feb2003} we study the question of existence and classification of reductions $H^{\flat}\to B$. 
We will see that the obstructions and the classification of these reductions up to equivalence fits completely with the classification {of}
extensions of $H\to E\to B$ to $T$-duality triples in terms of cohomological data {given in} \cite{MR2287642}.

\bigskip

A major part of the paper is devoted to the construction of a reduction $H^{\flat}\to E$ from a geometry on $H$. Geometric gerbes with band $\bbG$ have been introduced in \cite{MR2362847}, \cite{MR1876068}, or in \cite{MR1405064}, \cite{MR1794295} in the context of bundle gerbes.
In Section \ref{feb2010} we give a definition of a geometry on a gerbe using the language of stacks.
As a motivating example,  in Section \ref{feb1302} we make the induced geometry on a lifting gerbe explicit. This   essentially amounts to a translation of the results of \cite{MR1956150}, \cite{MR2029365}, \cite{MR2810948} from the bundle gerbe picture to the language used in the present paper. In Section \ref{feb2020} we give the construction of the reduction $H^{\flat}\to E$ from a geometry on $H\to E$. We show that every reduction can be obtained in this way, and we {analyse} under which conditions these reductions are fibrewise trivializable and hence good for $T$-duality.

\bigskip

In Section \ref{feb1920} we continue our example started in Section  \ref{feb1302}. We consider the canonical geometric lifting gerbe over a compact Lie group and use the explicit knowledge of its geometry in order to determine when the induced reduction is good for $T$-duality. The final result is stated as Theorem \ref{feb2560}.

\bigskip

This work was partially triggered by the recent paper \cite{2012arXiv1211.0763D} which observes that for a particular choice of the twist $H\to K\to K/T$ over a compact semisimple Lie group $K$ with only ADE-factors the $T$-dual is the underlying manifold of the  Langlands dual group $K^{L}$.  The authors of \cite{2012arXiv1211.0763D} uses a concept of $T$-duality based on differential forms.  In Section \ref{feb2040} we put this into the context of our integral $T$-duality and precisely describe the choice of the twist on $K$ which yields the Langlands dual group as the $T$-dual.

\bigskip

{\em Acknowledgement: We thank K. Waldorf guiding us through the  literature about the geometry on lifting gerbes.}

 \newcommand{\Obstr}{\mathrm{Obstr}}
 \newcommand{\Sh}{\mathrm{Sh}}
 \newcommand{\Red}{\mathbf{Red}}

\section{Classification of Reductions} \label{feb2003}

 In the present paper we use the language of stacks in smooth manifolds. 
 A good introduction is given in \cite{MR2206877}, see also \cite{MR2336247} for details about the sheaf theory in this framework. 

\bigskip

We consider the site $\Mf$ of smooth manifolds with the pretopology of open coverings.  
A manifold can {via the Yoneda embedding} be considered as a stack {on} $\Mf$. 
{An atlas of a stack $E$ on $\Mf$ is a map of stacks $M\to E$ from a smooth manifold $M$
which is a representable surjective submersion. A stack $E$ on $\Mf$ is called smooth if it admits an atlas.} 

\bigskip

The abelian category $\Sh_{\Ab}(B)$ of sheaves on a smooth stack $B$ is the abelian category of sheaves
on the site $\Mf_{sm/B}$  whose objects {$(M\to B)$} are representable submersions 
from manifolds to $B$ (cf. \cite[1.2.1]{MR2336247}). We will often omit the structure map to $B$ . If $\cF$ is a sheaf on $\Mf$, then we often use the symbol $\cF$ also for the restriction $\cF_{|\Mf_{sm/B}}$ in order to simplify the notation.

Let $E$ be a stack in $\Mf_{sm/B}$. The cohomology ${H^{*}(E;\cF)}$ of $E$ with coefficients in a sheaf $\cF\in \Sh_{\Ab}(B)$  
is defined in terms of the  derived functor of the evaluation functor  
$$ \Sh_{\Ab}(E)\to \Ab\ ,\quad \cF\mapsto \lim_{(M\to E)\in {\mathrm{Mor}}(\Mf_{sm/{B}})} \cF(M)\ .$$

\begin{rem2}
If $M\to E$ is an atlas, then we can form a Lie groupoid $M\times_{E}M\Rightarrow M$ and an associated simplicial manifold
$M^{\bullet}$ in $\Mf_{sm/{B}}$.  
 If $\cF $ is a sheaf on $B$, then we have a spectral sequence
\begin{equation}\label{apr2201}E_{1}^{p,q}\cong H^{p}(M^{q};\cF)\Rightarrow H^{p+q}(E;\cF)\end{equation}
converging to the cohomology of $E$ with coefficients in $\cF$.
{The differential $d_{1}:E_{1}^{p,q}\to E_{1}^{p,q+1}$ is derived from the simplicial structure of $M^{\bullet}$ in the usual manner.}

\bigskip

In the following we recall a few elements of sheaf theory on stacks from \cite{MR2336247}. Let $$i:\Sh_{\Ab}({B})\to \PSh_{\Ab}({B})$$ denote the embedding of sheaves into presheaves. A sheaf $\cF\in \Sh_{\Ab}({B})$ is called flabby if it is acyclic for the functor $i$.  If $\cF\in  \Sh_{{\Ab}}({B})$ is flabby, then in \eqref{apr2201} we have $H^{p}(M^{q};\cF)=0$ for all $p\ge 1$. Consequently the spectral sequence 
\eqref{apr2201} degenerates at $E_{2}$ for flabby sheaves. 
{
By  $\underline{\R}$ we denote the sheaf of rings of smooth {real-valued} functions {on} $\Mf$, i.e.  we have  $\underline{\R}(M)=C^{\infty}(M)$. If $\cF\in \Sh_{\Ab}({B})$ is a sheaf of $\underline{\R}$-modules, then it is flabby by the existence of smooth partitions of unity. }

\bigskip

{Let $\pi:E\to B$ be a representable {submersion}  between smooth stacks. For $(M\to B)\in \Mf_{sm/B}$ the composition  $M\times_{B}E\to E\to B$ is again an object of $\Mf_{sm/B}$.}
{We consider the functor $\pi_{*}:\Sh(B)\to \Sh(B)$ defined by $\pi_{*}(\cF)(M):=\cF(M\times_{B}E)$. If $\cF$ is  a} flabby sheaf on ${B}${, then it} is acyclic for the functor $\pi_{*}$ (\cite[Lemma 2.30]{MR2336247}).

\end{rem2}

\bigskip
{Throughout the paper the symbol $B$ will denote a proper smooth stack. Recall that a smooth stack $B$ is called proper if it admits an atlas $M\to B$ such that the Lie groupoid $M\times_{B}M\Rightarrow M$ is proper. Properness essentially says that the local stabilizer groups are compact.
\begin{rem2}
{
Smooth manifolds and 
 orbifolds are examples of proper smooth stacks. Another example is the stack $\cB K$ of $K$-principal bundles for a compact Lie group $K$. More general, if $G$ is a Lie group (not necessarily compact) which acts properly on a manifold $M$, then the quotient stack
$M/G$ is proper.}

{Our reason for assuming that $B$ is proper  is that {in this case} certain  sheaves on $B$ including $\underline{\R} $ are acyclic, see the proof of Lemma \ref{apr2202} for more details.}
\end{rem2}

\bigskip

{Assume that  $f:M\to B$ belongs to $\Mf_{sm/B}$}. Then the  tangent bundle of  $M$  has an integrable subbundle
$\cF_{{sm}/B}M\subseteq TM$. If $B$ is a smooth manifold, then  the fibre of $\cF_{sm/B}$  at $m\in M$ is defined by 
\begin{equation}\label{feb1320}\cF_{sm/B,m}M:=\{X\in T_mM\:|\:  df(X) =0\}\ .\end{equation}
 
 {
   In order to interpret \eqref{feb1320}  if $B$ is a smooth stack we proceed as follows. Let $ b:U\to B$ be an atlas of $B$. 
 Then we get a Lie groupoid $U\times_{B}U\Rightarrow U$ representing $B$. Its Lie algebroid is the vector bundle
 \begin{equation}\label{r34r34r34r34r3rerr34344r34r34r34r34}
\cL:=e^{*}\ker(d\pr_{0}:T(U\times_{B}U)\to TU) \ ,
\end{equation} where $e:U\to U\times_{B}U$ is the diagonal. The  differential  $	d\pr_{1}$ induces the anchor map, a morphism of Lie algebroids
 $a:\cL\to TU$. There exists a neighbourhood $V\subseteq M$ of $m\in M$ such that $f$ has a factorization $\tilde f$ over the atlas $b$.
   The technical definition of the fibre of the subbundle \eqref{feb1320} is the following:
 \begin{ddd}
 $$\cF_{sm/B,m}:=\{X\in T_mM\:|\:  d\tilde f(X) \in a(\cL_{\tilde \tilde f(m) }) \}\ .$$
 \end{ddd}
 One can check that this is independent of the choice of the atlas $U\to B$. If $B$ is smooth, then we can take $U=B$ and 
 revover  \eqref{feb1320}.}
 
 \bigskip
 
If $A$ is an abelian Lie group, then we can define the sheaf 
$\underline{A} $ on $\Mf$ of smooth functions to $A$. Its restriction to  $ \Mf_{sm/ B } $ contains the subsheaf $\underline{A}_{{sm}/B} $ of  sections 
  which are locally constant along $\cF_{sm}/B$.   
In detail {this subsheaf} is  given by  
\begin{equation}\label{feb2501}\underline{A}_{{sm}/B}(M ):=\{f\in  \underline{A}(M)\:|\: df_{|\cF_{{sm}/B}}=0\}\ . \end{equation}
 We write $A^{\delta}$ for the group $A$ considered as a discrete group.

\bigskip

We consider a diagram of maps between smooth stacks {on} $\Mf_{sm/{B}}$
\begin{equation}\label{feb2550}\xymatrix{H^{\flat}\ar[dr]\ar[rr]^{r}& &H \ar[dl]\\&E&}\ ,\end{equation}
where $H^{\flat}\to E$ is a  gerbe with band $\bbG_{{sm}/B}$, and $H\to E$ is a gerbe with band $\bbG $.
\begin{ddd}\label{feb1201}
We say that this diagram represents $H^{\flat}\to E$ as a $\bbG_{{sm}/B}$-reduction   of $H \to E$
if $r$ is a map of gerbes over $E$ over the canonical inclusion of  bands $\iota:\bbG_{{sm}/B}\to \bbG $.
 An equivalence between two $\bbG_{{sm}/B}$-reductions is a diagram     $$\xymatrix{ H^{\flat}\ar[d]^{r}\ar[rr]^{\tilde a} & & H^{\flat,\prime}\ar[d]^{r^{\prime}}\\H \ar[rr]^{a}\ar[dr]&  &H \ar[dl]\\&E&}\ ,$$
 where $a$ and $\tilde a$ are maps of gerbes with band $\bbG_{{sm}/B}$ and $\bbG $ over $E$, respectively. We let $$\Red(H\to E)$$ denote the set of equivalence classes of $\bbG_{{sm}/B}$-reductions of $H\to E$.
\end{ddd} 

The main goal of the present section is to calculate the set $\Red(H\to E)$ in cohomological terms. 
\bigskip

{Let $\cG\in \Sh_{\Ab}(B)$ be a sheaf of {abelian} groups on a smooth stack $B$ {and $\pi:E\to B$ be a map  of  stacks}. 
Then  the equivalence class of a gerbe $H\to E$ with band $\cG$
is classified by a characteristic class $[H\to E]\in H^{2}(E;\cG)$. 
In particular,
 the equivalence class of $H\to E$  as a gerbe over $E$ with band $ \bbG $ is determined by    the sheaf cohomology class
$[H\to E]\in H^{2}(E;\bbG )$. 

 Similarly, the equivalence class of  $H^{\flat}\to E$ as a gerbe over $E$ with band $\bbG_{{sm}/B}$ is determined by  the cohomology class $[H^\flat\to E]\in H^{2}(E;\bbG_{{sm}/B})$. 
The map $\iota: \bbG_{{sm}/B}\to \bbG $ of sheaves induces a map
$$\iota_*:H^{2}(E;\bbG_{{sm}/B})\to H^{2}(E;\bbG ) $$ of cohomology groups.
 The calculation of $\Red(H\to E)$ of based on the following observation.

\begin{lem}\label{feb1710}

The natural map $$\Red(H\to E) \to H^{2}(E;\bbG_{E/B})\ , \quad \eqref{feb2550}\mapsto [H^{\flat}\to E]$$ 
induces a bijection $$\Red(H\to E)\stackrel{\sim}{=} \iota_*^{-1}([H\to E])\ .$$
\end{lem}
\proof On the level of gerbes the map $\iota_*$ can be realized by a construction called  extension of bands. It  associates to a gerbe $H^{\flat}\to E$ with band $\bbG_{{sm}/B}$ the gerbe $ H^{\flat}\otimes_{B\bbG_{{sm}/B}}B\bbG \to E$ with band $\bbG $. The latter comes with a natural map
$$ i: H^{\flat}\to  H^{\flat}\otimes_{B\bbG_{sm/B}}B\bbG $$ 
 presenting $ H^{\flat}\to E$ as a 
$\bbG_{sm/B}$-reduction of $H^{\flat}\otimes_{B\bbG_{sm/B}}B\bbG \to E$. 

\bigskip

{In order to show the assertion we  construct an explicit inverse to the map $\Red(H\to E)\to \iota_*^{-1}([H\to E])$.
 {We consider {a} class $g \in \iota_*^{-1}([H\to E])$. If we  
 choose a $\bbG_{sm/B}$-gerbe $H^\flat\to E$ with $[H^{\flat}\to E]=g$, 
then we can find an equivalence  of $\bbG $-gerbes on $E$}
$$\alpha: H^{\flat}\otimes_{B\bbG_{sm/B}}B\bbG \stackrel{\sim}{\to}H\ .$$
We obtain a reduction  
\begin{equation}\label{may0901}
\raisebox{-4ex}{$(H^{\flat}\stackrel{\alpha\circ i}{\longrightarrow} H )\quad :=\quad $}
  \xymatrix{H^\flat\ar[dr]\ar[rr]^{\alpha \circ i}&& H\ar[dl]\\&E&}\ .\end{equation} 
We want to use this construction to define the desired inverse by
\begin{equation}\label{thomai5}
\iota_*^{-1}([H\to E]) \to \Red(H\to E) \qquad g \mapsto [H^{\flat}\stackrel{\alpha\circ i}{\longrightarrow} H ]\ .
\end{equation}
We have to show that this assignment is well-defined, in particular independent of the choice of $\alpha$. Therefore let $\alpha': H^{\flat}\otimes_{B\bbG_{sm/B}}B\bbG  \stackrel{\sim}{\to}H$ be another choice of equivalence. Then we get the diagram
$$
\xymatrix{
H^{\flat}\ar[d]^{i}\ar[rr]^{\id} & & H^{\flat}\ar[d]^{i}\\
H^{\flat}\otimes_{B\bbG_{sm/B}}B\bbG \ar[rr]^\id \ar[d]^\alpha & & H^{\flat}\otimes_{B\bbG_{sm/B}}B\bbG \ar[d]^{\alpha'}\\
H \ar[rr]^{\alpha'\circ \alpha^{-1}}\ar[dr]&  &H \ar[dl]
\\&E&}
$$
which shows that the reductions $(H^\flat \stackrel{\alpha \circ i}\longrightarrow H)$ and $(H^\flat \stackrel{\alpha' \circ i}\longrightarrow H)$ are equivalent. This shows that the assignment \eqref{thomai5} does not depend on the choice of $\alpha$. The functoriality of the extension of bands $H^\flat \mapsto H^{\flat}\otimes_{B\bbG_{sm/B}}B\bbG $ then implies that \eqref{thomai5} does also not depend on the choice of gerbe $H^\flat$. }

\bigskip

{
Thus the map \eqref{thomai5} is well-defined. Now we want see that it is indeed an inverse to the map $\Red(H\to E)\to \iota_*^{-1}([H\to E])$. Thus let $[H^\flat \stackrel{r}\to H]$ be in $\Red(H\to E)$. The morphism $r$ induces an equivalence 
$\alpha: H^{\flat}\otimes_{B\bbG_{sm/B}}B\bbG  \to H$. Then we have $[H^\flat \stackrel{r}\to H] = [H^\flat \stackrel{\alpha \circ i}\to H]$ in $\Red(H\to E)$. The remaining case is obvious. \hB
}

%
%
%
%
%
%
%
%

}


A pointed object of $\Mf_{sm/B}$ is a pair $(b,U )$  of an object
$U\to B$ in $\Mf_{sm/B}$ and a point $b\in U$. Let $\cF$ be a sheaf on ${B}$.
 For {a} cohomology class $x\in H^{*}(E;\cF)$  we consider its image \begin{equation}\label{feb1711}x_{(b,U\to B)} \in H^{*}( \pi^{-1}(b,U )  ;\cF):= \colim_{V}H^{*}(V\times_{B}E ;\cF) \ ,\end{equation}  where $V\subseteq U$ runs over the neighbourhoods of $b$. 
 We define $$F^{1} H^{*}(E;\cF):=\{x\in H^{*}(E;\cF)\:|\: (\forall\: \mbox{pointed objects $(b,U )$ in $\Mf_{|sm/B}$}\tho{\::\:} x_{(b,U )}=0)\}\ .$$

 We  define the sheaf
$\Omega^{*}_{sm/B}$ of complexes on $\Mf_{sm/M}$ whose evaluation on $M\to {B}$ is the complex 
$C^{\infty}(M,\Lambda^{*}\cF_{sm/B}^{*}M)$ of leafwise forms. 
By $\Omega^{1}_{sm/B,cl}\subseteq \Omega^{1}_{sm/B}$ we denote the subsheaf of closed {forms}. 

Let $\pi:E\to B$ be a representable submersion.
In order to calculate the cohomology $H^{*}(E;\Omega^{1}_{sm/B,cl})$  we use 
  the  resolution 
$$\Omega^{1}_{sm/B}\to \Omega^{2}_{sm/B}\to \dots\ .$$
\begin{lem}\label{apr2202}
For every $k\in \nat$ the sheaf $\Omega^{k}_{sm/B}\in \Sh_{\Ab}(B)$ is acyclic.
\end{lem}
\proof
Since $\Omega^{k}_{sm/B}$ is a sheaf of $\underline{\R} $-modules it is flabby.   Since $\pi$ is {a representable submersion}  the sheaf $\Omega^{k}_{sm/B}$ is acyclic for $\pi_{*}$.
From the Leray-Serre spectral sequence for $\pi$ we get an isomorphism 
$$H^{*}(E;\Omega^{k}_{sm/B})\cong H^{*}(B;\pi_{*}\Omega^{k}_{sm/B})\ .$$
Note that $\pi_{*}\Omega^{k}_{sm/B}$ is again a sheaf of $\underline{\R} $-modules and hence flabby. Let $U\to B$ be an atlas such that $U\times_{B}U\Rightarrow  U$ is a proper Lie groupoid. Then using the fact that the spectral sequence \eqref{apr2201} degenerates at $E_{2}$ we get an isomorphism
$$H^{*}(B;\pi_{*}\Omega^{k}_{sm/B})\cong H^{*}(\pi_{*}\Omega^{k}_{sm/B}(U^{\bullet}),d_{1})\ .$$
We get 
$$H^{p}(\pi_{*}\Omega^{k}_{sm/B}(U^{\bullet}),d_{1})=0\ , \quad \forall p\ge 1$$
from  Lemma \ref{apr2301}  applied to the following situation.
\begin{itemize}
\item The groupoid $\cG$ is the groupoid $U\times_{B}U\to U$.
\item The manifold $M$ is $U\times_{B}E$.
\item The total space of the vector bundle $V$ is
$U\times_{B} \Lambda^{k}({\cL\oplus} T^{v}\pi)^{{*}}$, where $T^{v}\pi$ is the vertical tangent bundle {of $M\to U$ and $\cL $ is the Lie algebroid, see  \eqref{r34r34r34r34r3rerr34344r34r34r34r34}.}
\item The complex $(C^{{\infty}}(M^\bullet,V^{\bullet}),d)$ is then identified with  $(\pi_{*}\Omega^{k}_{sm/B}(U^{\bullet}),d_{1})$.
\end{itemize} \hB 
Note that the full complex $\Omega^{*}_{sm/B}$ is an acyclic resolution of $\underline{\R}_{sm/B}$.
In particular we get for all $p\ge 1$ that
\begin{equation}\label{feb2551}H^{p}(E;\Omega^{1}_{sm/B,cl})\cong H^{p+1}(E;\underline{\R}_{sm/B} )\ .\end{equation}

\bigskip

We now consider the exact sequence of sheaves on {$B$}:
$$0 \to  \bbG_{sm/B}\to \bbG \stackrel{d\log}{\to} \Omega^{1}_{sm/B,cl}\to 0\ .$$
This together with \eqref{feb2551} gives a long  exact sequence in cohomology
\begin{equation}\label{feb1202}H^{1}(E;\bbG )\to H^{2}(E;\underline{\R}_{sm/B})\to H^{2}(E;\bbG_{sm/B})
\to H^{2}(E;\bbG ) \stackrel{\alpha}{\to} H^{3}(E;\underline{\R}_{sm/B})\ . \end{equation}

We define the subgroup
$$ H^{2}(E;\bbG )_{0} := \big\{ x\in H^{2}(E;\bbG )\:|\: \alpha(x)\in F^{1} H^{3}(E;\underline{\R}_{sm/B})\big\}\ .$$ 
 
The exponential sequence  induces an exact sequence of sheaves {on $B$}
\begin{equation}\label{apr2210}0\to \underline{\Z} \to \underline{\R}_{sm/B}\to \bbG_{sm/B}\to 0\ .\end{equation} The first map in this sequence induces the map in the denominator on the right-hand side of the following definition of the abelian group $A(E/B)$:
 $$A(E/B):=\frac{H^{2}(E;\underline{\R}_{sm/B})}{\im(H^{2}(E;\underline{\Z} )\to H^{2}(E;\underline{\R}_{sm/B}))}\ .$$
\begin{theorem}\label{feb1203}
Let $E\to B$ be a  {representable} submersion between smooth stacks and $H\to E$ be a gerbe with band $\bbG $. Then the following assertions hold true:
\begin{enumerate}
\item
The set $\Red(H\to E)$ is not empty if and only if
$[H\to E]\in H^{2}(E;\bbG)_{0}$.
\item 
If $\Red(H\to M)$ is not empty, then it is a torsor over the group
$A(E/B)$
\end{enumerate}
\end{theorem}
\proof
{We calculate $H^{p+1}(E;\underline{\R}_{sm/B} )$ using the resolution $\underline{\R}_{sm/B}\to \Omega_{sm/B}^{*}$.  In the proof of Lemma \ref{apr2202} we have seen that $R^{p}\pi_*\Omega^{{k}}_{E/B}$  vanishes for all $p\ge 1$ and  {$k\in \nat$}, {and it  is an acyclic}  sheaf on $B$} for  $p= 0$. 
Using the Leray-Serre spectral sequence  for $\pi$ we get the isomorphism $$ H^{p+1}(E;\underline{\R}_{sm/B} )  \cong
R^{p+1}\pi_*\underline{\R}_{sm/B}(B)\ .$$ 
In particular, if we want to check that an element of $H^{p+1}(E;\underline{\R}_{sm/B} )$ vanishes it suffices to show that its stalks at the points of $b$ vanish. This reasoning will be used below.

\bigskip

Let now $x\in H^{2}(E;\bbG ) $. We have $\alpha(x)=0$ in $R^{3}\pi_{*}\R_{sm/B}(B)$ if and only if for all pointed objects $(b, U )$ in $\Mf_{sm/B}$
we have
$\alpha(x)_{(b,U )}=0$, i.e. if and only if $x\in H^{2}(E;\bbG )_{0}$. 
%
%
In view of \eqref{feb1202} and Lemma \ref{feb1710}
 this proves Assertion 1.

\bigskip

Assertion 2. immediately follows from \eqref{feb1202} and the commutativity of the diagram
$$\xymatrix{H^{2}(E;\underline{\Z} )\ar[r]&H^{2}(E,\underline{\R}_{sm/B})\\
H^{1}(E;\bbG )\ar[u]^{\cong}\ar[r]&H^{1}(E;\Omega^{1}_{sm/B;cl})\ar[u]^{\cong}}\ .$$
{The left vertical map is an isomorphism in view of the long exact sequence associated to
\begin{equation}\label{apr2230}0\to \underline{\Z} \to \underline{\R} \to \bbG  \to 0\end{equation}
since $\underline{\R} $ is acyclic (Lemma \ref{apr2202} for $k=0$.)}
\hB 
We now assume that 
$H^{3}(\pi^{-1}(b,U ),\underline{\Z} )$ defined in \eqref{feb1711} is torsion free for all pointed objects $(b,U )$ of $\Mf_{sm/B}$.  This assumption is satisfied in the case of particular interest, where $E\to B$ is a principal $T$-bundle for a torus $T\cong (S^{1})^{n}$.  Under this assumption  the map
$$H^{3}(\pi^{-1}(b,U );\underline{\Z} )\to H^{3}(\pi^{-1}(b,U );\underline{\R} )$$ is injective and we get the left isomorphism in the diagram 
$$\xymatrix{ H^{2}(E;\bbG )_{0}\ar[r]\ar[d]^{\cong}&H^{2}(E, \bbG )\ar[d]_{\beta}^{\cong}\\
F^{1}H^{3}(E;\underline{\Z} ) \ar[r]&H^{3}(E;\underline{\Z} ) }\ .$$

\bigskip

We let $$\Red(H\to E)_0\subseteq \Red(H\to E)$$ denote the subset of reductions such that the underlying gerbes $ H^{\flat}\to E$ are trivial along the fibres of $E\to B$. 
 Our next goal is to calculate the set $\Red(H\to E)_0$. The answer will be given in terms of the fine structure of the Leray-Serre spectral sequence for $H^{*}(E;\underline{\Z} )$. This spectral sequence $(E_r,d_r)$ converges to the graded group $E_\infty^{p,q}\cong \Gr^{p} H^{p+q}(E;\underline{\Z} )$ with respect to a decreasing filtration $F^{*} 
H^{*}(E;\underline{\Z} )$.

\bigskip
 

We first draw a consequence of the absence of torsion in the fibrewise cohomology.

\begin{kor}\label{feb2601}
If in addition to the assumption of Theorem \ref{feb1203} the groups $H^{3}(\pi^{-1}(b,U );\underline{\Z} )$   are torsion-free for all pointed objects $(b,U )$ of $\Mf_{sm/B}$,
then the set $\Red(H\to E)$ is not empty if and only if 
$$\beta([H\to E])\in F^{1} H^{3}(E;\underline{\Z} )$$
\end{kor}
We have a natural map
$$\sigma:F^{1}H^{3}(E;\underline{\Z} )\to E_\infty^{2,1} \subseteq H^{1}(B;R^{2}\pi_*\underline{\Z} )$$
such that $F^{2}H^{3}(E;\underline{\Z} )=\ker(\sigma)$.
We further define the group \begin{equation}\label{feb2602}A(E/B)_0:= R^{2}\pi_*\underline{\Z} (B)/H^{{2}}(E;\underline{\Z} )\ .\end{equation}
By the Leray-Serre spectral sequence it fits into the exact sequence
$$0\to    \im(d_3^{0,2}:E_2^{0,2}\to E_2^{3,0}) \to A(E/B)_{{0}} \to \im(d_2:E_2^{0,2}\to E_2^{2,1}) \to 0\ .$$     
 
\begin{theorem}\label{feb1901} 
In addition to the assumption of Theorem \ref{feb1203} we assume {that $\pi:E\to B$ is proper and that}  the groups 
$H^{3}(\pi^{-1}(b,U ),\underline{\Z} )$ are torsion-free for all pointed objects $(b,U )$ of $\Mf_{sm/B}$.
We  further assume that $\Red(H\to E)$ is not empty.
\begin{enumerate}
\item The  set $\Red(H\to E)_0$ is not empty if and only if  $\beta([H\to E])\in F^{2}H^{3}(E;\underline{\Z} )$.
\item If $\Red(H\to E)_0$ is not empty, then it is a torsor over a group
$A(E/B)_0$
\end{enumerate}
\end{theorem}
\proof
We consider the class $x:=[H\to E]\in H^{2}(E;\bbG )_0$. By Theorem \ref{feb1203} we can find a lift $\tilde x\in H^{2}(E;\bbG_{sm/B})$  classifying a reduction in $\Red(H\to E)$.  The family of fibre restrictions $x_{(b,U )}$ defined by \eqref{feb1711} for all pointed objects $(b,U )$ in $\Mf_{sm/B}$  gives rise to a section $\{\tilde x\}\in R^{2}\pi_*\bbG_{sm/B}(B)$. The reduction classified by $\tilde x$ belongs to $\Red(H\to E)_0$ if and only if $\{\tilde x\}=0$.

\bigskip

We now use that $E\to B$ is {representable and a proper  submersion}  and hence a locally trivial fibre bundle with compact fibres. {For $p\in \nat$ we}
  have a locally constant sheaf $R^{{p}}\pi_*\underline{\R^{\delta}} $. 
It gives rise to a vector bundle $V^{{p}}\to B$ with flat connection $\nabla^{V^{{p}}}$ (often called Gauss-Manin connection) such that $R^{{p}}\pi_*\underline{\R^{\delta}}  $ is the sheaf of parallel sections of $V^{{p}}$. If $\cV^{{p}}$ denotes the sheaf of smooth sections of $V^{{p}}$, then the map
$ \underline{\R}^{\delta}  \to  \underline{\R}_{sm/B}$ induces an isomorphism
\begin{equation}\label{kjdkjljlqwdd}\cV^{{p}}\stackrel{\sim}{\to} R^{{p}}\pi_* \underline{\R}_{sm/B} \ .\end{equation}
The  map $\Z\to \R$ induces maps of sheaves
$R^{p}\pi_*
\underline{\Z} \to R^{p} \pi_*\underline{\R^{\delta} }$ 
whose images in $\cV^{{p}}$ we denote by $\cV^{p}_\Z$.
 The exponential sequence \begin{equation}\label{frf33f34f34f34f34f2344}
0\to \underline{\Z} \to \underline{\R}_{sm/B}\to \bbG_{sm/B}\to 0
\end{equation}

induces a sequence  of sheaves
$$0 \to \cV^{2}/\cV^{2}_\Z\to R^{2}\pi_*\bbG_{sm/B}\to R^{3}\pi_*\underline{\Z} \ .$$

  By assumption the image of $
 \{\tilde x\}$ in $R^{3}\pi_*\underline{\Z}(B) \cong {\Gr^{0} H^{3}(E;\underline{\Z})}$ vanishes so that we can consider $
 \{\tilde x\}\in (\cV^{2}/\cV^{2}_\Z)(B)$. 
 We have an exact sequence
 $$0\to \cV^{2}_\Z(B)\to \cV^{2}(B)\to \cV^{2}/\cV_\Z(B)\stackrel{\tilde \sigma}{\to} H^{1}(B;R^{2}\pi_*\underline{\Z} )\ .$$
Now we observe 
that
 $$\tilde \sigma (\{\tilde x\})=\sigma(x)\ .$$
 Hence $  \{\tilde x\}$ admits a lift to a section $ \lambda\in \cV^{2}(B)\cong {R^{2}}\pi_* \underline{\R}_{sm/B}$ if and only if $\sigma(x)=0$.
 
 \bigskip
 
Let us now assume that $\sigma(x)=0$ and let  
  $[\lambda]\in A(E/B)$ be the class of such a lift $\lambda$.
Then we form the class $\tilde x^{\prime}:=\tilde x-[\lambda]\in  \Red(H\to E)$. By construction we have
$\{\tilde x^{\prime}\}=0$. Hence $\tilde x^{\prime}\in \Red(H\to E)_0$.
This shows that first assertion.

\bigskip

Assume now that $\tilde x\in \Red(H\to E)_0$ and 
$\lambda\in \cV^{2}(B)$. Then
$\tilde x+[\lambda]\in  \Red(H\to E)_0$ if and only if
$\sigma([\lambda])=\tilde \sigma(\lambda)=0$, i.e.
$\lambda\in \cV^{2}_\Z(B)$.
 We conclude that
$\Red(H\to E)_0$ is a torsor over the 
group
$A_0(E/B)\cong \cV^{2}_\Z(B)/H^{2}(E;\underline{\Z} )$.
 \hB

\section{Geometry on  gerbes}\label{feb2010}

We consider a gerbe $H\to E$ with band $\bbG$. The goal of this subsection is to describe the axioms for a geometry on $H$.  If $W\to E$ is some stack over $E$, then  in order to simplify the notation we introduce the abbreviation $$W_{E}^{i}:=\underbrace{W\times_{E}\dots\dots\times_{E}W}_{{i \text{ times}}}$$ for the $i$-fold product of $W$ with itself over $E$.  The subscript in the  notation
$$\pr_{i_1\dots i_k}:W^{n}_E\to W^{k}_E$$ for the projection indicates in the standard manner 
on which factors it projects.  

\bigskip

Note that $H$ can be considered as the total space of a   $\bbG$-principal bundle 
$$H\cong H\times_{H}H\to H\times_{E}H= H^{2}_{E}\ .$$ In this case  
 we use the symbol $L$ in order to denote this  {circle} bundle as well as its total space.  
 The composition  of morphisms   is encoded in isomorphisms
$$\theta:\pr^{*}_{12}L\otimes \pr^{*}_{23}L\to \pr_{13}L$$ over $H_{E}^{3}$   which satisfies the cocycle condition
$$\pr^{*}_{134}\theta\circ ( \pr^{*}_{123}\theta\otimes \pr_{34}^{*}\id_{L})=\pr^{*}_{124}\theta \circ (\pr_{12}^{*}\id_{L}\otimes \pr^{*}_{234}\theta)\ .$$

    \begin{ddd}\label{feb1301}
A geometry on a $\bbG$-gerbe $H\to E$ consists of
\begin{enumerate}
\item  a $\Omega^{1}$-torsor $H^{geom}\to H$,
\item  a connection $\nabla$ on the $\bbG$-principal bundle  $L^{geom}\to H^{geom,2}_{E}$ defined 
by 
 pull-back 
 $$\xymatrix{L^{geom}\ar[r]\ar[d]&L\ar[d]\\H^{geom,2}_{E}\ar[r]&H^{2}_{E}}\ ,$$
 
\item a map $\omega:H^{geom}\to \Omega^{2}$ called curving.
 \end{enumerate}
 These data must satisfy the following conditions:
 \begin{enumerate}
 \item The pull-back of $\theta$ to $H^{geom,3}_{E}$ is parallel with respect to the connection 
 induced by $\nabla$.
   \item We have    $\pr_{1}^{*}\omega-\pr_{2}^{*}\omega=R^{\nabla}$, where $R^{\nabla}\in \Omega^{2}(H^{geom,2}_{E})$ denotes the curvature of $\nabla$.
     \item For $t,t^{\prime}\in H^{geom}(M)$ mapping to the same object in $E(M)$ and $\beta,\beta^{\prime}\in \Omega^{1}(M)$ we have \begin{equation}\label{ulimar16}\nabla^{t+\beta,t^{\prime}+\beta^{\prime}}\cong\nabla^{t,t^{\prime}}+\beta-\beta^{\prime}\ ,\end{equation}
   where $\nabla^{t,t^{\prime}}$ is the connection induced on the $\bbG$-principal bundle $$L(t,t^{\prime}):=(t,t^\prime)^{*}L^{geom}\to M\ .$$
   \item For  $t\in H^{geom}(M)$ and $\beta\in \Omega^{1}(M)$ we require that
 $\omega(t+\beta)=\omega(t)+d\beta$.
 \end{enumerate}  
 A geometric gerbe is a gerbe with a geometry. \end{ddd}
 
Let $t\mapsto \bar t$ denote the projection $H^{geom}(M)\to H(M)$. If $\beta\in \Omega^{1}(M)$ and $t\in H^{geom}(M)$, then $t+\beta\in H^{geom}(M)$ comes with a distinguished isomorphism $\bar t\to \overline{ t+\beta}$.
Implicitly, these isomorphisms are used in order to interpret \eqref{ulimar16}.
For a pair of sections $t,t^{\prime}\in H^{geom}(M)$ and a morphism $u:\bar t\to \bar t^{\prime}$ we can define the form 
$\nabla^{t^{\prime},t}\log u\in \Omega^{1}(M)$.

\bigskip

The curvature of the geometry is the form $\lambda\in \Omega^{3}(E)$
which is locally given by $d\omega(t)$ for any section $t$ of $H$. 
 
\bigskip 
 
 \begin{rem2}\label{apr0901}
 
The notion of a geometry on a gerbe with band $\bbG$ has been discussed in \cite{MR2362847}, \cite{MR1876068}, and in \cite{2012arXiv1211.0763D} in the bundle gerbe language.
For brevity {we} will explain the comparison of our definition  with \cite{MR1876068}. 
We assume that $E$ is a manifold. In Hitchin's picture \cite{MR1876068} a geometric gerbe was presented
with respect to an open covering $\cU:=(U_{\alpha})_{\alpha\in I}$ of $E$ in terms of
\begin{enumerate}
\item a
collection of $\bbG$-principal
 bundles $(L_{\alpha,\beta},\nabla^{L_{\alpha,\beta}}\to U_{\alpha,\beta})_{(\alpha,\beta)\in I^{2}}$
 with connection on double intersections, 
 \item a collection of bundle isomorphisms
 $(\theta_{\alpha,\beta,\gamma}:L_{\alpha,\beta}\otimes L_{\beta,\gamma}\to L_{\alpha,\gamma})_{(\alpha,\beta,\gamma)\in I^{3}}$ on triple intersections which preserve connections and satisfy the cocycle condition $$\theta_{\alpha,\gamma,\sigma}\circ( \theta_{\alpha,\beta,\gamma}\otimes \id_{L_{\gamma,\sigma}})=\theta_{\alpha,\beta,\sigma}\circ (\id_{L_{\alpha,\beta}}\otimes\theta_{\beta,\gamma,\sigma})
$$ on quadruple intersections, 
 \item and a collection of two-forms $(\omega_{\alpha}\in \Omega^{2}(U_{\alpha}))_{\alpha\in I}$ such that on double intersections we have 
 $\omega_{\alpha}-\omega_{\beta}=R^{\nabla^{L(\alpha,\beta})}$.
 \end{enumerate}
 
 {
 This data can be derived from a geometry in the sense of Definition \ref{feb1301} as follows.
We   choose sections $t_{\alpha}\in H(U_{\alpha})$   for all $\alpha\in I$. 
 Since $\Omega^{1}$ is a fine sheaf we have $H^{1}(U_\alpha;\Omega^{1})=0$. Hence the $\Omega^{1}$-torsor $t_\alpha^*(H^{geom})\to U_\alpha$ is trivial and we can choose sections $\hat t_\alpha \in H^{geom}(U_\alpha)$ that refine the $t_\alpha$. 
We get $\bbG$-principal bundles $L_{\alpha,\beta}:=L(  t_{\alpha},  t_{\beta})$ with connection $\nabla^{\hat t_{\alpha},\hat t_{\beta}}$. We further have isomorphisms
$\theta_{\alpha,\beta,\gamma}:=(t_{\alpha},t_{\beta},t_{\gamma})^{*}\theta$.
Finally we set   {$\omega_{\alpha}:=\omega(\hat t_{\alpha})\in \Omega^{2}(U_{\alpha})$}.
}

\bigskip
\newcommand{\Grbc}{\mathcal{G}rb^\nabla}
\newcommand{\calg}{\mathcal{G}}
Conversely  one can associate a geometry in the sense of Definition \ref{feb1301} to a 
Hitchin gerbe, or more generally to a bundle gerbe with connection. In the construction we use the fact that there is a 2-stack $\Grbc$ which associates to every smooth manifold $M$ a  $2$-groupoid $\Grbc(M)$ of gerbes with connection on $M$ as e.g. established in \cite{MR2318389} or \cite{MR2362847}. We do not go into the details but  only note the following key facts: 

\begin{itemize}
\item
 For every 2-form $\omega \in \Omega^2(M)$ there is a topologically trivial gerbe called $\mathcal{I}_\omega \in \Grbc(M)$\item 
Every $\bbG$-principal bundle $(L,\nabla)$ with connection induces  a morphism $f_{(L,\nabla)}: \mathcal{I}_\omega \to \mathcal{I}_{\omega+ R^{\nabla}}$ in $\Grbc(M)$. In particular, since a differential form $\alpha \in \Omega^1(M)$ can be considered as a connection on the trivial $\bbG$-bundle,  it induces  an automorphism 
 $f_\alpha:=f_{(M\times S^{1},d+\alpha)} :\mathcal{I}_\omega \to \mathcal{I}_{\omega+ d\alpha}$.
\item
The assignment $(L,\nabla) \mapsto f_{(L,\nabla)}$  extends to an   equivalence between the category of $\bbG$-bundles with connection of fixed curvature  and connection preserving morphisms and the category of morphisms  between the corresponding gerbes with connection. 
\end{itemize}

We consider a gerbe with connection $\calg \in \Grbc(E)$ such that its underlying gerbe $\calg^{top}$ is identified with $H$. More precisely,
$H(M)$ is the groupoid given as follows:
\begin{itemize}  \item
 Objects are pairs $(f,\phi)$, where $f:M\to E$ is a smooth map and $\phi:f^{*}\calg^{top}\to \cI^{top}$ is an isomorphism {to}
 the trivial topological gerbe $\cI^{top}$ on $M$ 
 \item A morphism   $(f,\phi)\to (f^{\prime}, \phi^\prime)$ can only exist  if $f=f^{\prime}$, and in this case it is a $2$-morphism $\phi\Rightarrow \phi^{\prime}$.
\end{itemize}
Note that the transition bundle $L$ on $H^2_E$ is the bundle which corresponds to the automorphism $$\cI^{top} \xrightarrow{\pr_1^*\phi^{-1}} s^*\calg^{top} \xrightarrow{\pr_2^*\phi} \cI^{top}$$
{where $s$ denotes the projection morphism $H^2_E \to {E}$. }
We now define the stack $H^{geom}$ which associates to $M$
 the groupoid $H^{geom}(M)$ as follows:
\begin{itemize}
\item Objects are triples $(f,\omega, \varphi)$ with $f:M\to E$ a smooth map, $\omega \in \Omega^2(M)$ and $\varphi: f^{*}\calg \to \mathcal{I}_\omega$ a morphism in $\Grbc(M)$.
\item A morphism $(f,\omega, \varphi) \to (f^{\prime},\omega',\varphi')$ in $H^{geom}(M)$ can   only exist  if $f=f^{\prime}$, $\omega = \omega'$, and in this case it is a 2-morphism $\varphi \Rightarrow \varphi'$ in $\Grbc(M)$.
\end{itemize}
We now indicate how $H^{geom}$ carries all the structure of a geometry. The action of $\Omega^1$ on $H^{geom}$ is provided by post composition: $\alpha \cdot (f,\omega, \varphi) := (f,\omega + d\alpha, f_\alpha \circ \varphi)$. 
The morphism $H^{geom}\to H$ is the obvious forgetful map, and 
 the map $H^{geom} \to \Omega^2$ is given by the projection $(f,\omega,\phi)\mapsto \omega$.
 {Finally the pullback of the transition bundle to $H^{geom,2}_{E}$
  refines to a morphism
 $$\mathcal{I}_{\omega_1} \xrightarrow{\pr_1^*\varphi^{-1}} s^*\calg \xrightarrow{\pr_2^*\varphi} \mathcal{I}_{\omega_2}$$
 in $\Grbc(H^{geom,2}_{E})$ between $\mathcal{I}_{\omega_1}$ and $\mathcal{I}_{\omega_2}$. Thus it can be identified with a bundle with connection.} 
 One can then check that these data establish a geometry. 
\end{rem2} 

Our motivation to work with the Definition \ref{feb1301} of a geometry is that this structure arises naturally   in the example of a lifting gerbe discussed in Section \ref{feb1302} and allows the simple construction of the reduction in Section \ref{feb2020}.

\bigskip

\begin{rem2}
The stack $H\to E$ has a presentation as $H^{geom}/\Omega^{1}$. We think  of $H^{geom}$ as the {``}minimal atlas''  for $H$
on which the curving $\omega$ is defined. 
Note that a typical model for the $\Omega^{1}$-torsor $H^{geom}\to H$  the sheaf of connections on some $\bbG$-principal bundle over $H$. From a more structural point of view
one can consider $H^{geom}\to E$ as a torsor over the Picard stack $B\bbG_{\nabla}$ of $\bbG$-principal bundles with connections.
\end{rem2}
\bigskip

In order to compare gerbes with geometries we introduce the notion of a morphism between geometric gerbes. We first introduce some notation. Let
$$\xymatrix{  
H\ar[rr]^{f}\ar[dr]&&H^{\prime}\ar[dl]\\
&E&} $$
be a morphism of gerbes over $E$ with band $\bbG$. Then we get a $\bbG$-bundle
$$Q(f):=(f,\id_{H^{\prime}})^{*}L^{\prime}\to H\times_E H^{\prime}$$
together with isomorphisms  
$$ 
\mu:\pr_{12}^{*}Q (f)\otimes \pr_{23}^{*}L^{ \prime} \stackrel{\sim}{\to}  \pr_{13}L^{\prime}\ , \quad  \nu: \pr_{12}^{*}L \otimes \pr_{23}^{*}Q (f)\stackrel{\sim}{\to} \pr_{13}^{*}Q (f) $$ of  $\bbG$-principal bundles on $ H \times_E H^{ \prime} \times_E H^{\prime}$ and
$H\times_EH\times_E H^{\prime}$ which are compatible with the composition maps $\theta$ and $\theta^{\prime}$ in the obvious manner.
If the two gerbes have a geometry and $F:H^{geom}\to H^{geom,\prime}$ lifts $f$, then  
$$Q^{geom}(f):= (F,\id_{H^{geom,\prime}})^{*}L^{\prime,geom}\cong \pr_{H\times_EH^{\prime}}Q(f)\to H^{geom}\times_EH^{\prime,geom}$$
has an induced connection $\nabla^{f}$ such that $\mu$ is parallel.

\begin{ddd} \label{feb1310}
An equivalence   between geometric gerbes   $$(f,F):(H^{geom}\to H,\omega,\nabla)\to (H^{\prime,geom}\to H^{\prime},\omega^{\prime},\nabla^{\prime})$$
 consists of 
 a commutative diagram of maps of stacks 
$$\xymatrix{H^{geom}\ar[rr]^{F}\ar[d]&&H^{\prime,geom}\ar[d]\\
H\ar[rr]^{f}\ar[dr]&&H^{\prime}\ar[dl]\\
&E&}\ ,$$ where $f$ is a morphism of gerbes over $E$ with band $U(1)$, and $F$ is a morphism of $\Omega^{1}$-torsors  over $f$
 such that \begin{enumerate}
 \item   $\nu$ induces an   isomorphism of $\bbG$-principal bundles with connection $$ \pr_{12}^{*}L^{geom}\otimes \pr_{23}^{*}Q^{geom}(f)\stackrel{\sim}{\to} \pr_{13}^{*}Q^{geom}(f)$$
on $H^{geom}\times_EH^{geom}\times_EH^{geom,\prime}$,
\item one $H\times_{E}H^{\prime}$ we have the identity of forms
$$ \pr_{H^{\prime}}^{*}\omega^{\prime}-\pr_{H}^{*}\omega=R^{\nabla^{f}}\ .$$\end{enumerate}
  \end{ddd}
{

\begin{rem2}
Using the comparison to other pictures of geometric gerbes as sketched in Remark \ref{apr0901} one can check that our notion of a geometric gerbes is equivalent the other pictures \cite{MR2362847}, \cite{MR1876068} {and \cite{MR1405064}}. This shows in particular  that equivalence classes of geometric gerbes are classified by the differential cohomology group $\widehat{H\Z}^{3}(E)$. 
\end{rem2}
 }

\section{The induced geometry on the lifting gerbe}\label{feb1302}

The details of the definition of a geometric gerbe {given in} Section \ref{feb2010}  were motivated by the example of the geometry on lifting gerbes. Using a different language, the geometry on lifting gerbes has been {thoroughly} discussed e.g. \cite{MR1956150}, \cite{MR2810948}, \cite{MR2029365}. In the present section we formulate these results in the language 
introduced in Section \ref{feb2010} and introduce the notation necessary for calculations in later sections.

\bigskip
 
We consider a central extension of Lie groups {(see Remark \ref{may1601} for the infinite-dimensional case)}
$$0\to U(1)\to \widehat G\to G\to 0\ .$$  
Our typical example is  the  central extension of the loop group  of a compact Lie group which will be further discussed in Section \ref{feb1920}. \begin{rem2}\label{may1601}
 In the context of the present paper the technical way to deal with  infinite dimensional groups of mappings like $G:=C^{\infty}
(S^{1},K)$ for a finite dimensional Lie group $K$ is to consider it as a sheaf of groups $\underline{G}$ on $\Mf$ using the exponential law, i.e. $\underline{G}(M)=C^{\infty}(M\times S^{1},K)$. In order to simplify the presentation  we will not use this language in the following section.
 \end{rem2}

\bigskip

Let $P\to E$ be a $G$-principal bundle. It determines a lifting gerbe $H\to E$ with band $\bbG$.  
 For a smooth manifold $M$  the groupoid $H(M)$ is the groupoid of pairs $(\phi,\hat P\to \phi^{*}P)$ of a map $ \phi:M\to E$ and a $\widehat{G}$-reduction $\hat P\to \phi^{*}P$ of the $G$-principal bundle $ \phi^{*} P$. For $t\in H(M)$ we write
$\phi(t)$ for the underlying map $M\to E$, $ P(t)$ for $\phi(t)^{*}P$, and $\hat P(t)$  for the total space of the corresponding $\widehat{G}$-principal bundle. 


\bigskip

Our goal is to explain how  geometric structures on the bundle $P$ induce a geometry on the lifting gerbe  $H\to E$. We essentially translate the construction first given by \cite{MR2029365}, \cite{MR1956150}, \cite{MR2810948}   from the language of bundle gerbes into the present picture.

\bigskip

The central extension of Lie groups induces a central extension of Lie algebras 
$$0\to \R\stackrel{\iota}{\to} \widehat{\gaaa}\stackrel{\tau}{\to} \gaaa\to 0\ .$$
We let $\sigma:\widehat{\gaaa}\to \R$ be {a}  left-inverse of $\iota$. The choice of $\sigma$ gives rise to a decomposition of vector spaces 

$$\widehat{\gaaa}\cong \R\oplus \gaaa\ .$$
As an abbreviation we will write $a \oplus_\sigma X$ for the element in $\widehat\gaaa$ corresponding to $(a,X) \in \R \oplus \gaaa$ under this isomorphism. 
Note that the adjoint action of $\widehat{G}$ on $\widehat{\gaaa}$ factors over $G$.
Following \cite{MR2029365}, \cite[Definition 3.15]{MR1956150}, \cite{MR2810948} we make the following definitions
\begin{ddd}  \label{feb1760}
 A splitting is map $$L:  P\times \widehat{\gaaa}\to \R$$ which is linear in the second variable and satisfies
  $$L(p\gamma,\widehat{\Ad}(\gamma^{-1})X)=L( p,X)\ , \quad \forall \gamma\in G\ ,  \forall X\in \widehat{\gaaa}\ , \forall p\in P$$
and  \begin{equation}\label{feb1303}L(\tilde p,\iota(r))=r\ , \quad\forall  r\in \R\ , {\forall \tilde p\in P} \ .\end{equation}
\end{ddd}

\begin{theorem}[\cite{MR1956150}, \cite{MR2810948}, \cite{MR2029365}]\label{feb1763}
The choice of  {a section }$\sigma$, a connection $\alpha$ and a splitting $L$ naturally determines a geometry on $H\to E$.
\end{theorem}
\proof
We only explain the construction of the data of a geometry according to Definition \ref{feb1301}.
The verification of the conditions are {straightforward} differential geometric calculations which are very similar to those given in the references.

\begin{enumerate}
\item

We first construct the torsor $H^{geom}\to H$. 
Let $p:\widehat{P}\to P$ be a $\widehat{G}$-reduction of $P$. We say that a principal bundle connection $\hat \alpha\in 
\Omega^{1}(\widehat{P})$ extends $\alpha$ if it satisfies the relation $\tau(\hat \alpha)=p^{*}\alpha$.
For a manifold $M$ we define $H^{geom}(M)$ to be the groupoid of pairs $(t,\hat \alpha)$ of a $\widehat{G}$-reduction $\hat P(t)\to \phi(t)^{*}P$ of bundles over $M$ and a connection $\hat \alpha$ on $\hat P$ which extends $\phi^{*}\alpha$. The set of morphisms   
$$\Hom_{H^{geom}(M)}((t,\hat \alpha),(t^{\prime},\hat \alpha^{\prime}))\subseteq \Hom_{H(M)}(t,t^{\prime})$$ is the subset of connection preserving isomorphisms of reductions.
The $\Omega^{1}$-torsor structure is given by 
$$(t,\hat \alpha)+\beta:=(t,\hat\alpha+\iota(\pr_{M}^{*}\beta))\ .$$

\item

Next we construct the connection $\nabla$ on the  bundle $L\to H^{geom,2}_{E}$. It suffices to describe the connection $\nabla((t,\hat \alpha),(t^{\prime},\hat \alpha^{\prime}))$ associated to a pair of sections
$((t,\hat \alpha),(t^{\prime},\hat \alpha^{\prime}))\in H^{geom,2}_{E}(M)$.
We first  provide  a description of the bundle $L(t,t^{\prime})\to M$. We  consider 
 the product
$$\hat P(t)\otimes \hat P(t^{\prime}) :=\hat P(t)\times_{ P(t),U(1)} \hat P(t^{\prime})\ , $$ where $U(1)$ acts on $\hat P(t^{\prime})$ with the opposite action.   The natural map
$$\hat P(t)\otimes \hat P(t^{\prime})\to {P(t)}$$ has the structure of an $U(1)$-principal bundle.
The diagonal action of $\widehat{G}$ induces a free action of $G$.
 The quotient is the $U(1)$-principal bundle
\begin{equation}\label{feb1304}L(t,t^{\prime}):= \hat P(t)\otimes \hat P(t^{\prime})/G \to M\ .\end{equation}

  Let now $(t,\hat \alpha)\in H^{geom}(M)$. We define a form in $\theta(t)\in \Omega^{1}(\hat P(t))$ by 
$$\theta(t,\hat \alpha):=\sigma(\hat \alpha(t))\ .$$
This is a connection one-form for the $U(1)$-bundle $\hat P(t)\to P(t)$. Note that in view of \eqref{feb1303} we have for $\beta\in \Omega^{1}(M)$ that
$$\theta((t,\hat \alpha)+\beta)=\theta(t)+\pr_M^{*}\beta \ .$$
The descent of the difference
   $$ \pr_{\hat P(t)}^{*}\theta(t)-\pr_{\hat P(t^{\prime})}^{*}\theta(t^\prime)\in \Omega^{1}(\hat P(t)\times_{ P(t)}\hat P(t^{\prime}))$$ to a form 
   $\delta(t,t^{\prime})\in \Omega^{1}(\hat P(t)\otimes  \hat P(t^\prime))$ exists and
  is a connection one-form on the $U(1)$-bundle $\hat P(t)\otimes  \hat P(t^\prime)\to P(t)$.
  By a straightforward calculation one can show that it is  
  $\widehat{G}$-invariant and consequently further descends along the quotient \eqref{feb1304} to a connection one-form $\delta((t,\hat \alpha),(t^{\prime},\hat \alpha^{\prime}))\in \Omega^{1}(L(t,t^{\prime}))$.
We define $\nabla^{(t,\hat \alpha),(t^{\prime},\hat \alpha^{\prime})}$ to be the corresponding connection on
$L((t,\hat \alpha),(t^{\prime},\hat \alpha^{\prime}))$.  
  
\item

We finally provide the curving $\omega$. It again suffices to describe
the forms $\omega(t,\hat \alpha)$ for sections $(t,\hat \alpha)\in H^{geom}(M)$.  

The map $\sigma$ is not $G$-equivariant. This is reflected by the
  map
\begin{equation}\label{feb1750}Z_{\sigma}:G\times \gaaa\to \R\ , \quad Z_{\sigma}(\gamma,X):=  \sigma(\widehat{\Ad}(\gamma) (0\oplus_{\sigma}  X)) \ , \quad  \gamma\in G\ ,\: X\in \gaaa\ .\end{equation}

The datum 
  of a splitting $L$ is equivalent to a map  (called reduced splitting in \cite[Definition 5.3]{MR2029365}, \cite[Definition 3.12]{MR1956150})
\begin{equation}\label{feb1751}l_\sigma:P\times \gaaa\to \R\end{equation} which is linear in the second argument and
satisfies
\begin{equation}\label{feb1761}l_\sigma(p\gamma,\Ad(\gamma^{-1}) X)=l_\sigma(  p,X)-Z_{\sigma}(\gamma^{-1},X) \ , \quad \forall \gamma\in G\ ,\forall X\in \gaaa\ , \forall p\in P\end{equation}
The relation between $L$ and $l_{\sigma}$ is given by 
\begin{equation}\label{feb1752}l_{\sigma}(p,X)=L(p,0\oplus_{\sigma}  X)\ .\end{equation}

Let $$\Omega(\alpha):=d\alpha+[\alpha,\alpha]$$ denote the curvature form associated to a connection one-form $\alpha$.
It turns out that the form  \begin{equation}\label{jan2910}\sigma(\Omega(\hat \alpha))+\pr_{ P}^{*}l_{\sigma}(\Omega(\alpha))\in \Omega^{2}(\hat P(t))\ . \end{equation}
  is horizontal and $\widehat{G}$-invariant. 
  The curving $\omega(t,\hat \alpha)\in \Omega^{2}(M)$ is uniquely determined by the condition that
  \begin{equation}\label{feb1850}\pr_{M}^{*}\omega(t,\hat \alpha)=\sigma(\Omega(\hat \alpha))+\pr_{ P}^{*}l_{\sigma}(\Omega(\alpha))\in \Omega^{2}(\hat P(t))\ .\end{equation}
\end{enumerate}
\hB

 One can calculate the curvature $\lambda\in \Omega^{3}(E)$ by calculating the derivative $d\omega(t)\in \Omega^{3}(M)$ for sections $(t,\hat \alpha)\in H^{geom}(M)$. Let $d_1 l_\sigma\in \Omega^{1}(P)\times \gaaa^{*}  $ denote the partial derivative of $l_\sigma$ with respect to the first variable. We get
 for $p\in P$ and horizontal vectors $X ,Y ,Z \in T^{h}_{p}P$
\begin{eqnarray}\label{jan2711} \lefteqn{\pr_{M}^{*}\lambda(X ,Y ,Z )}&&\\&=&    (d_1l_{\sigma})( X  , \Omega(\alpha)(Y ,Z ))-  (d_1l_{\sigma})(Y ,\Omega(\alpha)(X ,Z ))+(d_1 l_{\sigma})(Z ,\Omega(\alpha)(X  ,Y ))  \nonumber \ .\end{eqnarray}
From general principles it is clear that $\bbG$-principal bundles with connection 
act   on geometric gerbes.  In the example of the lifting gerbe it is easy to make this action explicit.

 \begin{theorem}
A  $\bbG$-principal  bundle with connection $S\to M$  naturally
  induces a morphism $$(f,F): (H^{geom}\to H,\omega,\nabla)\to (H^{geom}\to H,\omega+R^{S},\nabla)$$ of geometric gerbes in the sense of Definition \ref{feb1310}.    \end{theorem}
\proof
\begin{enumerate}
\item The map of stacks $f:H\to H$ sends
$t:=(\phi,\hat P(t)\to \phi^{*}P)$ to $f(t):=(\phi,\hat P{(t)}\otimes_{U(1)} S\to \phi^{*}P)$.
On morphisms it acts in the obvious manner.
\item The lift $F$ sends $(t,\hat \alpha)$ to $(f(t),\hat \alpha+_{U(1)}\pr^{*}_S\iota(\vartheta^{S}))$,
where $\vartheta^{S}\in \Omega^{1}(S)$ is the connection one-form on $S$ and
$ \hat \alpha+_{U(1)}\pr^{*}_S$ is a shorthand for the descent of the form
$\pr_1^{*}\hat \alpha+\pr_S^{*}\vartheta^{S}$ along the projection $\hat P(t)\times_M S\to \hat P{(t)}\otimes_{U(1)}  S$. 
\end{enumerate}
\hB
\section{Reduction via geometry}\label{feb2020}

We consider a {representable} submersion  $\pi:E\to B$  {over a proper smooth stack $B$}  and a  gerbe $H\to E$ with band $\bbG$. In this section we explain how a geometry on $H$ determines a reduction as introduced in Definition \ref{feb1201}.

\bigskip
 
The restriction   of the de Rham complex {$\Omega^{*}$} to {$\Mf_{sm/B}$}  has a decreasing filtration $F^{*}\Omega^{*} $ which in the evaluation {on} an object   $(M\to B)\in \Mf_{sm/{B}}$ is given by 
$$\hspace{-0.3cm}F^{p}\Omega^{n} (M )=\left\{\omega\in \Omega^{*}(M)\:|\:  \left(\begin{array}{c}\forall  m\in M\: \forall (X_{1},\dots,X_{{n-p+1}})\in (\cF_{sm/B,m}M)^{\times n-p+1}  \:|\: \\\iota_{X_{1}}\dots \iota_{X_{n-p+1}}\omega(m)=0\end{array}\right)\right\}\ .$$
Here $\cF_{sm/B}M\subseteq TM$ is the integrable subbundle introduced in \eqref{feb1320}. Note that this filtration induces the Leray-Serre spectral sequence for $H^{*}(E;\underline{\R^{\delta}} )$. In particular, elements in $F^{1}\Omega^{n} (M )$ are forms
which become trivial when restricted to fibres.

\begin{rem2}\label{fjljeflwefewfwf}  One should not confuse $\Omega^{*} $ with the fibrewise de Rham complex $\Omega^{*}_{sm/B}$ introduced in Section \ref{feb2003}. The  two complexes are related by a restriction map $\Omega^{*} \to \Omega^*_{sm/B} $.  
\end{rem2}

 \begin{theorem}\label{feb0301}
We assume that the gerbe $H$ is a equipped with a geometry as in Definition \ref{feb1301} whose curvature satisfies $\lambda\in F^{1}\Omega^{3} (E)$. Then we can define a canonical 
reduction 
$$\left( \xymatrix{H^{\flat}_{{}_{} }\ar[dr]\ar[rr]^{a}&& H \ar[dl]\\&E&} \right)\in \Red(H \to E)\ .$$ 

\end{theorem}
\proof
The geometry provides the $\Omega^{1}$-torsor $H^{geom}\to H$, the curving $\omega\in \Omega^{2}(H^{geom})$, and the  bundles $L(t,t^{\prime})\to M$ with connections $\nabla(t,t^{\prime})$ for pairs of sections $t,t^{\prime}\in H(M)$ over the same object of $E(M)$.
 \begin{ddd}
Let $M\to B$ be  in $\Mf_{sm/B}$.\begin{enumerate}
\item 
A section $t\in H^{geom}(M)$ is called {fibrewise flat } if $\omega(t)\in F^{1}\Omega^{2} (M)$.
\item Let $t,t^{\prime}\in H^{geom} (M)$. A morphism $u:\bar t\to \bar t^{\prime}$ in $H(M)$ is called {fibrewise flat}, if  
 $$\nabla^{t^{\prime},t}\log u\in F^{1} \Omega^{1} (M )\ .$$ 
\end{enumerate}
\end{ddd}

 We now define a  stack  $H^{\flat} $  {on} $\Mf_{sm/B}$ {as follows}.
 Let $M\to B$  be an object of $\Mf_{sm/B}$. \begin{enumerate}
\item 
An object $t\in  H^{\flat}(M )$ is  a {fibrewise} flat object in $H^{geom}(M)$. 
\item A morphism $u:t\to t^{\prime}$ between {fibrewise} flat objects
is a {fibrewise} flat morphism $u:\bar t\to \bar t^{\prime}$.
The composition of morphisms in $H^{\flat}$ is induced from the composition of morphisms in $H$.
 \end{enumerate}
 One now checks that the prestack given by this description is in fact a stack {but this is not essential for the rest of the argument, since we could alternatively replace it by its stackification}.
By construction we have a diagram of stacks
 $$
 \xymatrix{&H^{geom}\ar[dr]&\\H^{\flat}\ar[dr]\ar[rr]&&H\ar[dl]\\&E&}\ .
$$
\begin{lem}
 $H^{\flat}\to E$ is a gerbe with band $\bbG_{sm/B}$.
\end{lem}
\proof
{If $M\to B$ is an object in $\Mf_{sm/B}$,  $f\in E(M)$,}   and  $m\in M$, then we must show that there is an open neighbourhood $V\subseteq M$ of $m$ such that the fibre   $ H^{\flat}(V)_{/f_{|V}} $ is not empty.
Since $H\to E$ is a gerbe 
we can find an open neighbourhood $V\subseteq M $ of $m$  and an object
$\bar t\in H(V)_{/f_{|V}} $. After shrinking $V$ if necessary there exists a lift $ t \in H^{geom}(V)$.
Then $f_{|V}^{*}\lambda  =d\omega(t)\in F^{1}\Omega^{3} (V ) $.   After shrinking $V$ further and using a foliated  Poincar\'e lemma we can find   a form
$\beta\in \Omega^{1}(V )$ such that $\omega(t )+d\beta\in F^{1} \Omega^{2} (V )$.    Then $$ t +\beta \in  H^{\flat}( V  )_{/f_{|V}} \ .$$
 
\bigskip

We now show that two sections of $H^{\flat}$ are locally isomorphic.
Assume that $t,t^{\prime}\in  {H}^{\flat }(M )_{{/f}} $ are two objects and $m\in M$.
Let $\bar t,\bar t^{\prime}\in H(M)_{{/f}}$ be the projections. Since $H\to E$ is a gerbe  there exists a neighbourhood $U$ of $m$ and an isomorphism $\bar u:\bar t_{|U}\to \bar t^{\prime}_{|U}$ in $H(U)_{/f_{|U}}$. We have
$$d\nabla^{t^{\prime},t}\log u=R^{\nabla^{t^{\prime},t}}=\omega( t^{\prime}_{|U} )-\omega( t_{|U})\in F^{1 }\Omega^{2 } (U )\ .$$
After shrinking $U$ if necessary we can find a section $\alpha\in \bbG (U)$ such that
 $$ \nabla^{t^{\prime}_{|U}, t_{|U}}  \log u +d\log \alpha\in F^{1}\Omega^{1 } (U )\ .$$  
 It follows that $\alpha u:\bar t_{|U}\to \bar t^{\prime}_{|U}$ is a flat morphism from $t_{|U}$ to $t^{\prime}_{|U}$. 
 
\bigskip 
We finally calculate the band of $H^{\flat}\to E$. We fix $t\in H^{\flat}(M )$. Then we have a canonical identification
$\Aut_{H(M)}(\bar t)_{{/\id_{f}}}\cong \bbG(M)$. An automorphism $u\in \Aut_{H(M)}(\bar t)_{{/\id_{f}}}$ is flat if and only  {if}
$u\in \bbG_{sm/B}(M )\subseteq \bbG(M)$.   Hence $\Aut_{H^{\flat}(M )}(t)_{{/\id_{f}}}\cong \bbG_{sm/B}(M )$.
\hB

In Theorem \ref{feb1203} we calculated the set $\Red(H \to E)$. It is now a natural question how this calculation is related with the construction of reductions via geometry  given in \ref{feb0301}. The precise answer is given by the following proposition.

\begin{theorem}\label{feb2460}
We assume that $\pi:E\to B$ is a proper  submersion of smooth stacks and that
$H\to E$ is a gerbe with band $\bbG$.
\begin{enumerate}
\item 
The set $\Red(H \to E)$ is not empty if and only if there exists a geometry on $H$ with curvature $\lambda\in F^{1}\Omega^{3} (E)$.
\item\label{pointtwo} Every element of  $\Red(H \to E)$ can be realized as the 
  reduction determined by some geometry according to Theorem  \ref{feb0301}.
\end{enumerate}
\end{theorem}
\proof
If such a geometry exists, then we have constructed a reduction   in Theorem \ref{feb0301}.

\bigskip

Conversely, assume that a reduction exists. By  Theorem \ref{feb1203}  we have $ [H \to E] \in  H^{3}(E;\bbG )_{0}$. We choose some geometry on the gerbe $H\to E$ with curvature $\lambda\in \Omega^{3}(E)$. Let $\bar \lambda\in \Omega^{3}_{sm/B}(E) $  denote the restriction {(see Remark \ref{fjljeflwefewfwf})} of $\lambda$ to $\Lambda^{3}\cF_{sm/B}E$. 
We have (with $\alpha$ as in \eqref{feb1202})
$$0=\alpha ([H \to E] )_{(b,U )}=[\bar \lambda]_{(b,U )}\in  H^{3}(\pi^{-1}(b,U );\underline{\R}_{sm/B})$$
for all pointed objects $(b,U )$ in $\Mf_{sm/B}$.
As we have seen before in the proof of Theorem \ref{feb1203} this implies that
$[\bar \lambda]=0\in H^{3}( B;\underline{\R}_{sm/B})$.
Hence there exists a form $\bar\alpha\in \Omega^{2}_{sm/B}(E)$ with $d\bar \alpha=\bar \lambda$. Let $\alpha\in \Omega^{2}(E)$ be any lift of $\bar \alpha$. 
Then $$\lambda-d\alpha \in F^{1}\Omega^{3} (E)   \ .$$  We now change the geometry of $H$
by adding $-\alpha$ to the curving. This provides a new geometry such that $\lambda\in F^{1}\Omega^{3} (E)$. {This finishes the proof of 1.}

\bigskip
  
Assume now that we have found one reduction \begin{equation}\label{feb1604}
 {x=\left(\xymatrix{H^{\flat}_{{}_{} }\ar[dr]\ar[rr]&& H \ar[dl]\\&E&}\right)\in \Red(H \to E)\ .}
\end{equation} associated to some geometry. By Theorem \ref{feb1203}
any other can be obtained from this by shifting with a section $\gamma\in  R^{2}\pi_*\underline{\R}_{sm/B}(B)$. One can find a form $\alpha\in \Omega^{2}(E)$ whose restriction $\bar \alpha\in \Omega^{2}_{sm/B}(E)$ is closed and  represents $\gamma$. Then we can shift the geometry by adding $\alpha$ to the curving. The associated reduction is  the shift $x+\gamma\in \Red(H \to E)$.              
This can be seen by an argument which is similar to that for Lemma \ref{feb0501} below.
\hB

Let $H\to E$ be a gerbe with band $\bbG$ which is equipped with a geometry with curving $\lambda\in F^{1}\Omega^{3} (E)$. According to Theorem \ref{feb0301}
we get a reduction of the form \eqref{feb1604}. In the study of $T$-duality the following condition plays an important role:
\begin{ddd}\label{feb1720}We say that $H^{\flat}\to E$ is fibrewise trivializable, if 
  for each pointed object $(b,U )$ of $\Mf_{sm/B}$ there exists a neighbourhood $V\subseteq U$ of $b$ such that  the gerbe
$V\times_{B}H^{\flat}\to V\times_{B}E$ is trivializable.
\end{ddd}
 
 \begin{rem2}
 {The name `fibrewise trivializable' is a bit misleading here since it suggests that the restriction $H^{\flat}_{|E_{b}}$ 
 is a trivializable gerbe with band $\underline{S}^{1,\delta}$  for every point $b:*\to B$ of $B$. 
 But this condition does not even make sense, since in general the inclusion $\iota_{b}:E_{b}\to E$ of the fibre of $E$ over $b$  {is not a morphism in the site  $\Mf_{sm/B}$ since $E_{b}\to B$ is not a submersion.}
 The restriction $H^{\flat}_{|E_{b}}\to E_b$ could be understood as a gerbe over the band $\iota_{b}^{*}\bbG_{sm/B}$, but this would lead to a different notion and is 
 not what we want. }
 \end{rem2}

 The following theorem gives a geometric characterization of  
 condition of fibrewise trivializability.


\begin{theorem}\label{feb1770}
We adopt the assumptions of Theorem \ref{feb0301}.
We further assume that   $B$ is connected and that $V\times_{B} H^{\flat}\to V\times_{B}E$  is trivializable for some   object $V\to B$ of $\Mf_{sm/B}$. Then $H^{\flat}\to E$ is fibrewise trivializable
  if and only    if  {$\lambda\in F^{2}\Omega^{3}(E)$}.  \end{theorem}
\proof
Let $(b,U )$ be a pointed object of $\Mf_{sm/B}$. 
The gerbe 
  $U\times_{B}H^{\flat}\to U\times_{B}E$ is classified by a class
$$\kappa:=[U\times_{B}H^{\flat}\to U\times_{B}E]\in H^{2}(U\times_{B}E;\bbG_{E/B})\cong b{R^{2}\pi_{*}\underline{S}^{1}_{sm/B}(U)} \ .$$
The long exact sequence for $R^{*}\pi_{*}$ associated to \eqref{apr2230} and the fact the
$\underline{\R} $ is flabby and hence $\pi_{*}$-acyclic provides  isomorphisms
$$R^{p}\pi_{*}\bbG \cong R^{p+1}\pi_{*}\underline{\Z} $$ for all $p\ge 1$.
From  the short exact sequence
$$0\to \bbG_{sm/B}\to \bbG \to \Omega^{1}_{sm/B,cl}\to 0$$ 
we now get the exact sequence {(see \eqref{kjdkjljlqwdd} for $\cV^{2}$)}
 \begin{equation}\label{feb1701}R^{2}\pi_*\underline{\Z}  \stackrel{\gamma}{\to}  \cV^{2}\stackrel{a}{\to}   R^{2}\pi_*\bbG_{sm/B}\stackrel{\beta}{\to}   R^{3}\pi_*\underline{\Z} \ .\end{equation}

The sheaf $$\cV^{2}\cong R^{1}\pi_*(\Omega^{1}_{sm/B,cl})\cong H^{2}(\pi_* \Omega_{sm/B})$$ 
(here we again use that $\Omega^{1}_{E/B,cl}\to \Omega^{\ge 1}_{E/B}[1])$ is a flabby, and hence $\pi_{*}$-acyclic resolution)
is the sheaf of sections of a  real vector bundle with a  natural flat connection $\nabla$, the Gauss-Manin connection. 

\bigskip

  Note that  $R^{3}\pi_*\underline{\Z} $ is locally constant. We consider  a section $\kappa\in  
 R^{2}\pi_*\bbG_{sm/B}(U)$. If   the stalk of $\beta(\kappa)\in R^{3}\pi_*\underline{\Z} (U)$ vanishes in one point $b\in U$, then $\beta(\kappa)$ vanishes in a whole neighbourhood. We replace $U$ by such a neighbourhood. In addition we will assume that $U$ is contractible.
Then $\kappa=a(\psi)$ for some section $\psi\in \cV^{2}(U)$. We can consider its derivative
$$\nabla \psi\in (\Omega^{1} \otimes_{C^{\infty}} \cV^{2})(U)\ .$$
Observe that  the image of $\gamma$ consists of parallel sections. 
If $\psi$ is parallel and $a(\psi)(b)=0$ for one point $b\in U$, then $a(\psi)=0$.
Conversely, if $a(\psi)=0$, then $\psi$ is parallel.
 
 \bigskip

{In the rest of the proof we interpret the fact that $\psi$ is parallel in terms of the condition on the curvature form  $\lambda$  of the gerbe.}
 We assume that the section $\kappa \in R^{2}\pi_*\bbG_{sm/B}( {U})$ corresponds to a reduction $x$ as in \eqref{feb1604}.  Let $V\subseteq U$ be a neighbourhood of $b$ such that the gerbe $V\times_{B}H^{\flat}\to V\times_{B}E$ is trivializable and choose  a section $t_V\in H^{\flat}(V\times_{B}E) $ {with image $\bar t\in H(V\times_{B}E)$}. Since $\beta(\kappa)=0$ and $U$ is contractible we  can assume that $\bar t_V$ {extends   to a section $\bar t\in H (U\times_{B}E )$. We can furthermore find a lift $t\in H^{geom}(U\times_{B}E )_{/\bar t}$,} not necessarily flat. If $\bar\omega(t)\in \Omega^{2}_{sm/B}(U\times_{B}E)$ denotes the restriction of $\omega(t)\in \Omega^{2}(U\times_{B}E)$ to the vertical foliation $\cF_{sm/B}(U\times_BE)$, then  we have $d\bar \omega(t)=0$ by our assumption on the curvature $\lambda$ of the geometry. 
 Hence we get a class $[\bar\omega(t) ]\in H^{2}(U\times_{B}E  ;\underline{\R}_{sm/B})$.
    \begin{lem}\label{feb0501}
We have
$\kappa =a([\bar \omega(t)]) \in H^{2}(U\times_{B}E ;\underline{{S}}^{1}_{sm/B})$.
\end{lem}
\proof
We can choose a sufficiently fine open covering $\cU:=(U_{i})_{i\in I}$ of $U\times_{B}E$ so that the following argument works.
 For 
 each $i\in I$ there exists a form $\beta_{i}\in \Omega^{1}(U_{i})$ which satisfies $d\bar \beta_{i}=\bar \omega(t)_{| U_{i}}$, where again $\bar \beta\in \Omega_{sm/B}^{1}( U_{i})$ is the restriction of $\beta$ to the vertical foliation.
   Then $t_{|  U_{i}}-\beta_{i}\in H^{\flat}( U_{i} )_{{/t_{|U_{i}}}}$.  
   
   \bigskip

Note that $ ( \bar \beta_{i}-\bar \beta_{j})_{| U_{i}\cap U_{j}}$ is closed. Hence we  can
 choose a collection of sections $\alpha_{ij}\in \bbG( U_{i}\cap U_{j})$ indexed by pairs  $(i,j)\in I\times I$ such that   $$ d\log\alpha_{ij}- (\beta_{i}-\beta_{j})_{| U_{i}\cap U_{j}}\in F^{1}\Omega^{1} (U_{i}\cap U_{j} )\ .$$   
 Then  $$
\nabla^{t+\beta_{j},t+\beta_{i}} \log(\alpha_{ij}\id_{\bar t_{|U_{i}\cap U_{j}}})  \in F^{1}\Omega^{1} (U_{i}\cap U_{j} )
  \ ,
$$
In other words,     $f_{ij}:=\alpha_{ij}\id_{\bar t_{|U_{i}\cap U_{j}}}$ is a  morphism from $(t+\beta_{i})_{|U_{i}\cap U_{j}}$ to $(t+\beta_{j})_{|U_{i}\cap U_{j}}$  {in}
$H^{\flat}( U_{i}\cap U_{j} )_{{  /t_{|U_{i}\cap U_{j} }}}$ over $\id_{\bar t_{|U_{i}\cap U_{j}}}$
for each  pair $(i,j)\in I^{2}$.

\bigskip

We form the \v{C}ech cocycle $(c_{ijk})\in C^{2}(\cU,\bbG_{sm/B})$ by \begin{eqnarray*} 
 c_{ijk} &{:=}&f_{ij}\circ f_{jk}\circ f_{ki}\\&=&(\alpha_{ij}\alpha_{jk}\alpha_{ki}\id_{{\bar t}},0)\in \Aut_{H^{\flat}(U_{i}\cap U_{j}\cap U_{k} \to E)}((t+\beta_{i})_{|U_{i}\cap U_{j}\cap U_{k}})  \\&\cong& \bbG_{sm/B}( U_{i}\cap U_{j}\cap U_{k} ) \ .\end{eqnarray*}
 By definition, the cocycle   $(c_{ijk}) $    represents the class $\kappa \in H^{2}(U\times_{B}E ,\bbG_{sm/B})$.
On the other hand, this class is by construction the image under $a$ of the   class $[\bar \omega(t) ]\in H^{2}(U\times_{B}E,\underline{\R}_{sm/B})$. \hB

By Lemma \ref{feb0501} we can take 
$$\psi:=[\bar \omega(t)] \in H^{2}( \Omega^{*}_{sm/B}(U\times_{B}E))\cong R^{2}\pi_{*}\underline{\R}_{sm/B}(U )\ .$$

We now use  a connection on the bundle  $U\times_{B}E\to U$ in order to define horizontal lifts $X^{h}\in \Gamma(E_{b^{\prime}},T(U\times_{B}E))$ of tangent vectors $X\in T_{b^{\prime}}U$. The Gauss-Manin connection is given by 
$$\nabla_X\psi =[\overline{\cL_{X^{h}}\omega(t)} ] =[\overline{\iota_{X^{h}} d+d \iota_{X^{h}} \omega(t)} ]=[\overline{\iota_{X^{h}} d\omega(t)}]\ , \quad X\in \Gamma(U,TU)$$
and independent of the choice of the connection above.
From $d\omega(t)=\lambda$ we get
$$\nabla_X\psi=[\overline{\iota_X^{h}\lambda}] \ .$$
 
We have $\lambda\in F^{2}\Omega^{3} (U\times_{B}E)$ if and only if $\nabla\psi$ is parallel.
 \hB 





\section{Canonical construction of $T$-duality triples}\label{feb2002}

Let $T$ be a torus, i.e. an abelian Lie group isomorphic to $(S^{1})^{n}$ for some $n\in \nat$.
We consider a principal $T$-bundle $E\to B$ over some smooth stack $B$. In particular, $E\to B$ is a proper submersion. Furthermore, we consider  a gerbe $H^{\flat}\to E$ with band $\bbG_{sm/B}$, see \eqref{feb2501}.
\begin{ddd}\label{feb1903}  On the site $\Mf_{sm/B}$ we define the stack  $\hat H^{{\flat}}$  {put a ${}\flat$ in order to save $\hat H$ for the extension to $\bbG$.} of fibrewise trivializations of $H^{{\flat}}\to E$ such that its evaluation on $(M\to B)\in \Mf_{sm/B}$  is the groupoid 
$$\hat H^{\flat}(M ):=H^{\flat}(M\times_{B}E )\ .$$
\end{ddd}
Since $M\to B$ {and $ E\to B$ are submersions, the induced map 
$M\times_{B}E\to E\to B$ is a submersion, too.}  Therefore  the stack $H^{\flat}$ can be evaluated on this object.

 \bigskip

We consider the factorization \begin{equation}\label{feb1610}\hat H^{\flat}\stackrel{\hat q}{\longrightarrow} \hat E\stackrel{\hat \pi}{\longrightarrow} B\ ,\end{equation} where $\hat E$ is the relative moduli stack of $\hat H^{\flat}\to B$. 
In the special case when $B$ is a manifold, the stack  $\hat E$ is simply the sheaf of sets {on} $\Mf_{sm/B}$ associated to the presheaf 
$\pi_{0}\hat H^{\flat}$ of isomorphism classes in $\hat H^{\flat}$.
 The main goal of the present section is to describe the structure of the sequence \eqref{feb1610}.

\bigskip

We consider the dual torus $$\hat T:=\Hom(\pi_1(T,1),S^{1})$$
of $T$  as a Lie group. 

\begin{theorem}
We assume that the gerbe $H^{\flat}\to E$  is  fibrewise trivial{izable} (Definition \ref{feb1720}). Then
the projection $\hat \pi:\hat E\to B$ is naturally  a principal $\hat T$-bundle and $\hat q:\hat H^{\flat}\to \hat E$ has the structure of a gerbe with band  {$\bbG_{\hat E/B}$}. 
\end{theorem}
\proof
We consider a pointed object $(b,M )$ in $\Mf_{sm/B}$.
By Definition \ref{feb1720} there exists a neighbourhood $U\subseteq M$ of $b$ such that  the gerbe $U\times_{B}H^{\flat}\to U\times_{B}E$ is trivializable. {We consider $\tilde t=(U\times_{B}E\to E)\in E(U\times_{B}E)$. It has a lift
$t\in H^{\flat}(U\times_{B}E)_{/\tilde t}$ which represents its isomorphism class $\bar t\in \hat E(U)$. 
This shows that $\hat E$ is not empty locally over $B$.}

 \begin{lem}
 There exists a natural identification of sheaves of abelian groups
\begin{equation}\label{feb2510} \nu:R^{1}\pi_{*} \bbG_{sm/B}\cong \underline{\hat T} \ .\end{equation}
\end{lem}
\proof {
We use the exact sequence \eqref{frf33f34f34f34f34f2344}. Since the integral cohomology of $T$ is torsion-free we get a short exact sequence
$$0\to R^{1}\pi_{*} \Z \to R^{1}\pi_{*} \underline{\R}_{sm/B}\stackrel{p}{\to} R^{1}\pi_{*} \bbG_{sm/B}\to 0\ .$$ In order to calculate $R^{1}\pi_{*} \underline{\R}_{sm/B}$ we use the resolution
 (see Lemma \ref{apr2202}) $$\Omega^{0}_{sm/B} \to \Omega^{1}_{sm/B}\to\Omega^{2}_{sm/B}\to \dots\ .$$} 
We consider an object $(M\to B)\in \Mf_{sm/B}$ such that
the bundle $M\times_{B}E\to M$ has a section $s:M\to M\times_{B}E$.
  Let $\psi\in R^{1}\pi_{*} \R_{sm/B}(M )$ be represented by $\omega\in \Omega^{1}_{sm/B,cl}(M\times_{B}E )$. Then we define
  $$\nu({p}(\psi))\in \underline{\hat T}(M )\ , \quad \nu({p}(\psi))(m)([\tau])= [\int_{{S^{1}}} \tau(m)^{*}\omega]\in S^{1} \ .$$
  Here $m\in M$, $\nu(\psi)(m)\in \Hom(\pi_1(T,1),S^{1})$, 
    and the path $\tau(m):  S^{1}\to M\times_{B}E$ is given by  $\tau(m)(t):=s(m)\tau(t)$,  where $\tau:S^{1}\to T$ represents $[\tau]\in \pi_1(T,1)$.   One checks that $\nu({p}(\psi))(m)$ is well-defined independent{ly} of the choice representatives {$\psi$ of $p(\psi)$,}  $\tau$ of $[\tau]$  and the section $s$. \hB 
  
 Since $H^{\flat}\to E$ is a gerbe with band $\bbG_{sm/B}$ the set of isomorphisms classes in the fibre $\hat H^{\flat}(M )\to \hat E(M)$ is  empty or 
 naturally a torsor over the group $H^{1}({M\times_{B}E};\bbG_{sm/B})$. 
 The sheafification of the presheaf $$(M\to B)\mapsto H^{1}({M\times_{B}E };\bbG_{sm/B})$$ is
 naturally isomorphic to  
 $$R\pi_{*}^{1}\bbG_{sm/B} \stackrel{\eqref{feb2510}}{\cong} \underline{\hat T}  \ .$$ 
{Since $\hat E$ is locally not empty over $B$} we see that  $\hat E\to B$ is a principal $\hat T$-bundle.
  
 \bigskip


Let $t\in \hat H^{\flat }(M )$ {and $\bar t\in \hat E(M)$ be its image. }
Since $H^{\flat}\to E$ is  a gerbe with band $\bbG_{sm/B}$ and the fibres of $M\times_{B}E\to M$ are connected we have 
$$\Aut_{\hat H^{\flat}(M)} (t)_{/\bar t}\cong \bbG_{sm/B}(M\times_{B}E )\cong    {\bbG_{sm/B}(M)}\ .$$
This shows that $\hat q:\hat H^{\flat}\to \hat E$ is a {again a} gerbe with band  {$\bbG_{ sm/B}$}. \hB 

%

%

 
Over $E\times_{B}\hat H $ we have the tautological section $\taut$. Hence we get a diagram  of stacks {on} $\Mf_{sm/B}$:
 \begin{equation}\label{feb1730}\xymatrix{&&H^{\flat}\times_{B}\hat H^{\flat} \ar[dr]^{\hat r}\ar[dl]^{r}&&\\&H^{\flat}\times_{B}\hat E\ar[dr]\ar[dl]&&E\times_{B}\hat H^{\flat}  \ar[dr]^{\hat q}\ar@{.>}@/_1cm/[ul]^{\taut} \ar[dl]&\\ H^{\flat}\ar[dr]^{q}&&E\times_{B} \hat E\ar[dl]^{p}\ar[dr]^{\hat p}&&\hat H^{\flat} \ar[dl]^{\hat q}\\&E\ar[dr]^{\pi}&&\hat E\ar[dl]^{\hat \pi}&\\&&B&&}\ .\end{equation}
   We call \eqref{feb1730} the universal trivialization diagram.

 {We let $H\to  E$ and $\hat H\to \hat E$ denote the gerbes obtained by extension of the band along $\bbG_{sm/B}\to \bbG$, see the proof of Lemma \ref{feb1710}. The tautological section induces a morphism of gerbes with band $\bbG$}
%
{
  $$\xymatrix{H\times_{B}\hat E\ar[dr]&&\ar[ll]^{u} E\times_{B}\hat H \ar[dl]\\&E\times_{B}\hat E&}\ .$$}
We get an induced diagram{\scriptsize  \begin{equation}\label{feb1740}\hspace{-1cm}\xymatrix{
& H \times_{B}\hat E\ar[dr]\ar[dl] &&E\times_{B}\hat H  \ar[dr]^{\hat q}\ar[ll]^{u}  \ar[dl]&\\
H \ar[dr]^{q}&&E\times_{B}\hat E\ar[dl]^{p}\ar[dr]^{\hat p}&&\hat H \ar[dl]^{\hat q}\\
&E\ar[dr]^{\pi}&&\hat E\ar[dl]^{\hat \pi}&\\
&&B&&}\ .\end{equation}}
 \begin{ddd}\label{feb1902}
We call \eqref{feb1740} the canonical $T$-duality diagram associated to the pair $H^{\flat}\to E\to B$.
\end{ddd} 
 \bigskip

  \bigskip
 
 \begin{rem2} 
One can check that the  canonical $T$-duality diagram is the restriction
to $\Mf_{sm/B}$ of a $T$-duality diagram in the sense of \cite[Definition 4.1.3]{MR2797285}.
If one interprets gerbes with band $\bbG$ as twists, then by   \cite[Lemma 4.1.5]{MR2797285} equivalence classes
$T$-duality diagrams over a manifold $B$ are in bijection with
 $T$-duality triples introduced in  \cite[Definition 2.8]{MR2287642}.
If $B$ is a space, then the classification of {$\bbG_{sm/B}$-} reductions of a gerbe $H\to E$ with band $\bbG $ over a torus bundle $E\to B$
given in Theorem \ref{feb1203} fits completely with the classification of $T$-duality triples over $B$  given in \cite[Theorem 2.24]{MR2287642}.  
\end{rem2}
\bigskip

\bigskip 
 
 \begin{rem2}
 {In this remark we show how the Chern classes of the dual torus bundle $\hat E$ depend on the choice of the reduction $H^{\flat}$ of $H$.}
We fix an isomorphism $T\cong (S^{1})^{n}$ and a generator $y\in H^{1}(S^{1};\Z)$. In this way we fix
classes $y_{i}\in H^{1}(T;\Z)$ for $i=1,\dots,n$ such that $H^{*}(T;\Z)\cong \Z[\lambda_{1},\dots,\lambda_{n}]$.

\bigskip

 Let $E\to B$ be a principal $T$-bundle over a connected base and $(E_{r}^{p,q},d_{r})$ be the associated Leray-Serre spectral sequence. Then we can identify $E_{2}^{0,1}\cong H^{1}(T;\Z)$ and therefore view the classes $y_{i}$ as elements of $E_{2}^{0,1}$. The Chern classes of $E\to B$ are defined by
$$c_{i}(E\to B):=d_{2}(y_{i})\in E^{2,0}_{2}\cong H^{2}(B;\Z)\ , \quad i=1,\dots,n\ .$$
Our choices also fix dual generators $\hat y\in H^{1}(\hat T;\Z)$ and therefore the Chern classes $\hat c_{i}(\hat E\to B)\in H^{2}(B;\Z)$, $i=1,\dots,n$ of a principal $\hat T$-bundle $\hat E\to B$.

\bigskip


Let us fix a gerbe $H\to E$ with band $\bbG$. 
The underlying gerbe with band $\bbG_{sm/B}$ of a reduction $x\in \Red(H \to E)$ will be denoted by $H^{\flat}(x)\to E$.    We assume that $x\in \Red(H \to E)_{0}$ 
and consider a representative $v\in \cV^{2}_{\Z}(B)$ 
of a class $[v]\in A_{0}(E/B)$ in the group defined in \eqref{feb2602}.    By Theorem \ref{feb1901} we can form the reduction $x^{\prime}:=x+[v]\in  \Red(H \to E)_{0}$.
Finally we let $\hat E\to B$ and $\hat E^{\prime}\to B$ denote the underlying principal $\hat T$-bundles
obtained as $T$-duals associated to the reductions $H^{\flat}(x)\to E$ and $H^{\flat}(x+[v])\to E$. 

\bigskip

We define a collection of integers $(B_{ij})_{1\le i<j\le n}$ such that
$$v=\sum_{1\le i<j\le n}B_{ij}\:y_{i}\wedge y_{j}\ .$$

\begin{prop}
 We have
 $$\hat c_{k}(\hat E^{\prime}\to B)-\hat c_{k}(\hat E\to B)= \sum_{1\le i<k} B_{ik}c_{i}(E\to B)-\sum_{k<j\le n}B_{kj}c_{j}(E\to B)\ . $$ 

 \end{prop}
 This follows from 
\cite[Theorem 2.24, 2.]{MR2287642} \hB 
\end{rem2}

\section{The transgression gerbe of the tautological $K$-bundle}\label{feb1920}

In this section we give an illustration of the general constructions of the present paper by working out an explicit example.

\bigskip

Let $K$ be a compact connected  semisimple Lie group with Lie algebra $\kaaa$. By $LK$ we denote its loop group, i.e. the group  $C^{\infty}(S^{1},K)$
of smooth maps from $S^{1}$ to $K$ {(see Remark \ref{may1601})}. We  fix a central extension
\begin{equation}\label{feb2301}1\to S^{1}\to \widehat{LK}\to LK\to 1\ .\end{equation}
The plan of this section is as follows:

\begin{enumerate}
\item
We construct a tautological $LK$-principal bundle $\tilde P\to K$ and obtain a lifting gerbe
$H\to K$ with band $\bbG$ by the construction of Section \ref{feb1302}.

\item
We show that
the bundle $\tilde P\to K$ has a canonical connection and splitting. We therefore obtain a geometry on the lifting gerbe $H\to K$ by Theorem \ref{feb1763} 

\item

Let $T\subseteq K$ be a maximal torus. We consider $K$ as the total space of a principal $T$-bundle  $K\to B$ over the flag manifold $B:=K/T$.  We use the  geometry on the lifting gerbe $H\to K$ in order to construct a reduction $H^{\flat}\to K$ by Theorem \ref{feb0301}.

\item

We finally investigate  under which conditions   $H^{\flat}\to K$ is fibrewise trivializable in the sense of Definition  \ref{feb1720} so that  we can form a $T$-duality triple {using} the construction of Section \ref{feb2002}. 
\end{enumerate}

In order to state 
our final result   Theorem \ref{feb2560} we need 
the following construction. 
 Since $T$ is abelian  the exponential map 
  $\exp:\taaa\to T$ is a homomorphism of groups. Its kernel
   $\Lambda:=\ker(\exp)\subset \taaa$ is called the integral lattice.
   We use  the notion  $t\mapsto [t]$ for the projection
   $\R\to  \R/\Z \cong S^{1}$ {and} define 
   the {group} homomorphism 
\begin{equation}\label{apr1120}\Lambda\ni \lambda\mapsto  \gamma_{\lambda}\in LK\ , \quad \gamma_{\lambda}([t]):=\exp(t\lambda)\ , \quad t\in \R\ .\end{equation}
We then  define a central extension
\begin{equation}\label{feb2430}0\to S^{1}\to \hat \Lambda\to \Lambda\to 0\end{equation}
of $\Lambda$ as a pull-back of the central extension \eqref{feb2301} along the homomorphism \eqref{apr1120}.

\begin{theorem}\label{feb2560}
For each central extension \eqref{feb2301} of the loop group $LK$ of a connected semisimple compact Lie group $K$ {there is an associated  geometric lifting gerbe
$H\to K$. The construction of $H \to K$ is {sketched}   in the points 1.-3. above and made precise in the proof below. Its curvature is the unique bi-invariant three form  in the respective cohomology class in $H^3(K,\R)$}.

{The geometry of $H \to K$ induces a reduction which is fibrewise trivial along the principal $T$-bundle $K\to K/T$ if and only  if 
the  central extension \eqref{feb2430} is trivial. }
  \end{theorem}
\proof

\bigskip
\centerline{1.}
\bigskip

We start with the construction of $\tilde P\to K$ by a specialization of a general construction called   transgression. Let
  $$\pi:P\to E\times S^{1}$$ be a  principal $K$-bundle over some stack $E$ such  that $P$  is trivializable along the fibres of  the projection $E\times S^{1}\to E$.  Under this assumption 
 we can define a  principal $LK$-bundle $$\tilde \pi:\tilde P\to E\ ,$$
called the transgression of $P\to E\times S^{1}$.  By definition  the evaluation $\tilde P(M)$ of the stack $\tilde P$  on the object $M\in \Mf$  is the groupoid   $P(M\times S^{1})_{/\id_{S^{1}}}$. 

\bigskip

Therefore, in order to construct the principal $LK$-bundle $\tilde P\to K$ we must construct a principal $K$-bundle $P\to K\times S^{1}$ which is fibrewise trivializable. To this end
we first consider the trivial bundle
 $$K\times K\times \R\to K\times \R\ , \quad (h,g,t)\mapsto (g,t)\ .$$  The structure of a  principal $K$-bundle is given by the right action of the group $K$ on the first factor. The group $\Z$ acts on this bundle by 
  $$n\cdot(h,g,t):=(g^{n}h,g,t+n)\ ,\quad n\cdot (g,t):=(g,t+n)\ , \quad n\in \Z\ .$$
 We define the principal $K$-bundle  $P\to K\times S^{1}$ by taking the quotient with respect to this $\Z$-action.
 Since $K$ is connected the bundle $P\to K\times S^{1}$ is trivializable along the fibres of the projection $K\times S^{1}\to K$.

 \bigskip
 
 Let $\tilde P\to K$ be the transgression of $P$ and $H\to K$ be the associated lifting gerbe with band $\bbG$.
 
 \begin{rem2}
A point $\tilde p\in \tilde P$ in the fibre over $g\in K$ is a twisted loop
  in $K$, i.e.
a map
$\tilde p:\R\to K$ such that $\tilde p(t+n)= g^{-n}\tilde p(t)$ for all $t\in \R$ and $n\in \Z$. 
\end{rem2}
 
\bigskip

\centerline{2.}

\bigskip

Our next task is to describe an essentially canonical connection and a splitting for $\tilde P\to K$ in the sense of  Definition \ref{feb1760}. It will be derived from a connection $\alpha$ on $P\to K\times S^{1}$ by {the following}  construction.

\bigskip

We again start with the more general case of the transgression $\tilde \pi:\tilde P\to E$ of a bundle $\pi:P\to E\times S^{1}$. In addition, 
we assume that we are given  a connection on $P$ represented by a connection one-form $$\alpha\in \Omega^{1}(P)\otimes\kaaa\ .$$ 
We then obtain an induced connection on $\tilde P\to E$ with {a} connection one form $$\tilde \alpha\in \Omega^{1}(\tilde P)\otimes L\kaaa\ .$$ In order to describe $\tilde \alpha $, using an atlas of $E$,  we can assume that $E$ is a manifold. Let $b\in E$ and $\tilde P_{b}\cong \Gamma(\{b\}\times S^{1},P)$ be the fibre of $\tilde P$ over $b$. Note that for a section $\tilde p\in \tilde P_b$ we have an identification $$T_{\tilde p}\tilde P\cong \left\{\tilde X\in \Gamma(S^{1},\tilde p^{*}TP)\:|\: \Big(\exists X\in T_bE\:\: \forall t\in S^{1}\:|\: d(\pr_E\circ\pi)(\tilde X(t))=X\Big)\right\}\ .$$  With this notation $d\tilde \pi(\tilde X)=X$. The form $\tilde \alpha$ is given by 
$$\tilde \alpha(\tilde X)(t)=\alpha(\tilde X(t))\ , \quad t\in S^{1}\ ,$$
where we adopt the identification $L\kaaa\cong C^{\infty}(S^{1},\kaaa)$.

\bigskip

Following \cite{MR1956150} {and} \cite{MR2810948} we now describe the splitting.
First of all, the connection $\alpha$  gives a map
$$A_{\alpha}:\tilde P\to L\kaaa \ , \quad A_{\alpha}(\tilde p):= (\tilde p^{*} \alpha)(\partial_{t})$$
which is  sometimes called the Higgs field, see \cite[Definition 6.5]{MR1956150}, \cite[Equation (15)]{MR2029365}.

We fix a {vector space} split $\sigma:\widehat{L\kaaa}\to \R$ of the extension of Lie algebras
\begin{equation}\label{feb2402}0\to \R\to \widehat{L\kaaa}\stackrel{a}{\to} L\kaaa\to 0\ .\end{equation}
The structure of central extensions of loop groups of compact Lie groups is well-known \cite{MR900587}. In particular, by   \cite[Proposition 4.2.4]{MR900587}
we can  {assume}   that the   cocycle 
$$K_{\sigma}\in \Lambda^{2}L\kaaa^{*}\ , \quad K_{\sigma}(X,Y):=\sigma([0\oplus_{\sigma} X,0\oplus_{\sigma} Y])\ , \quad X,Y\in L\kaaa$$
is of the form $$K_\sigma(X,Y)=\int_{S^{1}} \langle dX,Y\rangle$$
for a uniquely determined 
symmetric invariant bilinear form \begin{equation}\label{feb2440}\langle\dots ,\dots\rangle :\kaaa\otimes \kaaa \to \R\ .\end{equation}
 This compares well with \cite[Equation (37)]{MR1956150}.
The following Lemma  is taken from \cite[Propositions 6.1 and 6.2]{MR1956150}.
\begin{lem}[Gomi]
The connection $\tilde \alpha$ provides a canonical reduced splitting (see \eqref{feb1751}) by
\begin{equation}\label{feb2210}l_{\sigma}(\tilde p,X):=  \int_{S^{1}} \langle A_{\alpha}(\tilde p),X\rangle \ dt\ .\end{equation}
\end{lem}

We let $L$ denote the splitting given by $l_{\sigma}$ and \eqref{feb1752}.

\bigskip

Therefore, in order to define a connection $\tilde \alpha$ and a splitting for $\tilde P\to K$ we must define a connection $\alpha$ on $P\to K\times S^{1}$. We first define a $\Z$-invariant connection on the principal $K$-bundle
 $K\times K\times \R\to K\times \R$.
 We choose a function $\chi\in C^{\infty}([0,1])$ which is
constant near the endpoints of the interval  and satisfies $\chi(0)=0$ and $\chi(1)=1$.
Recall that we denote points in $K\times K\times \R$ by $(h,g,t)$. We use the suggestive notation $dgg^{-1}, h^{-1}dh\in \Omega^{1}(K)\otimes \kaaa$ for the left- and right invariant {Maurer-Cartan} forms
and $h Xh^{-1}$ for $\Ad(h)(X)$, $h\in K$, $X\in \kaaa$.
Then we let $$\alpha\in \Omega^{1}(K\times K\times [0,1])\otimes \kaaa$$ be given by
$$\alpha:h^{-1}dh-\chi(t) h^{-1} dg g^{-1}h\ .$$
This is a connection one-form on the bundle
$K\times K\times [0,1]\to K\times [0,1]$ which can be extended to a 
$\Z$-invariant connection one-form on
$K\times K\times \R\to K\times \R$.
In this way we get the connection $\alpha$ on the quotient $P\to K\times S^{1}$.

\bigskip

By the general construction explained above we get an essentially canonical connection $\tilde \alpha$ and a splitting $L$ for $\tilde P\to K$. The only choice is that of the cut-off function $\chi$ which turns out to be   irrelevant in the following. We get an induced geometry
on the lifting gerbe $H\to K$ by Theorem \ref{feb1763}.

\bigskip

\centerline{3.}

\bigskip

We now consider the $T$-principal bundle $K\to K/T=B$. We want to construct a reduction of the gerbe $H \to K$ with band $\bbG$ to a gerbe $H^{\flat}\to K$ with band $\bbG_{sm/B}$ using Theorem \ref{feb0301}. To this end we must check that the curvature  {$\lambda$} of the geometry of $H\to K$ satisfies
$\lambda\in F^{1}\Omega^{3} (K )$. 

\bigskip

Hence we must calculate the curvature of the geometric gerbe $H\to K$. The invariant symmetric bilinear form $\langle\dots ,\dots \rangle$ on $\kaaa$ induces an element 
$$c\in (\Lambda^{3}\kaaa^{*})^{K}\ , \quad c(X,Y,Z):=\langle [X,Y],Z\rangle\ , \quad X,Y,Z\in \kaaa\ .$$  
\begin{lem}\label{uliapr0401}
The curvature $\lambda\in \Omega^{3}(K)$ of the canonical geometry on the lifting gerbe $H\to K$ is the bi-invariant three-form  on $K$ determined by $c$.
\end{lem}
\proof
Let $\Omega(\alpha)\in \Omega^{2}(P)\otimes \kaaa$ denote the curvature form of $\alpha$.
By a direct calculation, the horizontal derivative of $A_\alpha$ evaluated at $\tilde X\in T_{\tilde p}\tilde P$ is given by 
\begin{equation}\label{jan2710}(d^{h}A_{\alpha})(\tilde p)(\tilde X)(t)= \Omega(\alpha)(\tilde X(t),d\tilde p(t)(\partial_{t})) \ , \quad t\in S^{1} \ .\end{equation}
Using \eqref{jan2710}, \eqref{feb2210} and \eqref{jan2711} we can calculate the curvature of the gerbe.
Let $X_{i}\in \kaaa$ for $i=1,2,3$,  $k\in K$ and let $kX_{i}\in T_{k}K$ be the corresponding tangent vectors. We write $(k,t)X_{i}, (k,t)\partial_{t}\in T_{(k,t)}(K\times S^{1})$ for the corresponding tangent vector at $K\times S^{1}$ and  abbreviate $\Omega:=\Omega(\alpha)$.
\begin{equation}\label{jan3001}\lambda(kX_{1},kX_{2},kX_{3})=\frac12 \sum_{\sigma\in \Sigma_{3}} \sign(\sigma)
\int_{S^{1}}\langle \Omega((k,t)X_{i_{\sigma(1)}},(k,t)\partial_{t}), \Omega((k,t)X_{i_{\sigma(2)}},(k,t)X_{i_{\sigma(3)}})\rangle  .\end{equation}
By a direct calculation we get 
$$\Omega =-d\chi\wedge \Ad(h^{-1}) dg g^{-1}+(\chi^{2}-\chi)\Ad(h^{-1})[ dgg^{-1},dgg^{-1} ]\ .$$
This gives
\begin{eqnarray*}
\int_0^{1}\langle \Omega((k,t)X_{i_{1}},(k,t)\partial_{t}), \Omega((k,t)X_{i_{2}},(k,T)X_{i_{3}})\rangle&=&2
 \int_{0}^{1} \chi^{\prime}(t)(\chi^{2}(t)-\chi(t))dt \:\langle X_{i_{1}},[X_{i_{2}},X_{i_{3}}]\rangle \\
&=&-\frac{1}{3} c(X_{i_{1}},X_{i_{2}},X_{i_{3}})
\end{eqnarray*}
This formula implies the assertion.
\hB 
\begin{kor}\label{feb2211}
We have $\lambda\in F^{1}\Omega^{3} (K )$.
\end{kor} 
\proof
The bundle $\pi:K\to K/T$ admits a left action of $K$. 
Note that $\lambda$ is left invariant. The torus $T$ is a particular fibre of the bundle $K\to B$.
If $\lambda_{|T}=0$, then $\lambda_{|kT}=0$ for all $k\in K$ and hence $\lambda \in F^{1}\Omega^{3} (K )$. Since $\taaa$ is abelian, we have $c_{|\Lambda^{3}\taaa}=0$. This implies $\lambda_{|T}=0$. 
\hB 

By Corollary \ref{feb2211} and Theorem \ref{feb0301}
we can  define a reduction 
$$\left(\xymatrix{H^{\flat}\ar[rr]\ar[dr]&&H \ar[dl]\\
&K&}\right)\in \Red(H \to K)$$

This finishes the proof of the first part of the Theorem.

\bigskip

\centerline{4.}

\bigskip

For the second assertion of the Theorem 
 we {analyse} when the reduction $H^{\flat}\to K$ constructed in  step 3. is fibrewise trivializable. First we check the condition on the curvature assumed in Theorem \ref{feb1770}.

\begin{lem}
We have $\lambda\in F^{2}\Omega^{3} (K )$.
\end{lem}
\proof
We again use that the bundle $\pi:K\to B$ is left invariant and that $\lambda$ is invariant.
Hence we must check, that $\iota_{X}\iota_{Y}\lambda=0$ for $X,Y\in T_{e}T$. This is equivalent to 
  the condition that $c(X,Y,Z)=0$ for $X,Y\in \taaa$ and $Z\in \kaaa$.
But this is clear since $c(X,Y,Z)=\langle [X,Y],Z\rangle$ and $\taaa$ is abelian.
\hB

Next we reduce the question whether $H^{\flat}\to K$ is fibrewise trivializable to the question 
on the maximal torus.

 \begin{lem}\label{apr1101}
If there exists a section $t\in H^{geom}(T)$ with $\omega(t)=0$, then
there exists a neighbourhood $V\subseteq B$ of $[T]\in B$ such that $H^{\flat}_{|\pi^{-1}(V)}\to \pi^{-1}(V)$ is trivial.
\end{lem}
\proof
We choose a contractible neighbourhood $V\subseteq B$  of $[T]$. There exists a section $\check{t}\in H^{geom}(\pi^{-1}(V))$ such that $\check{t}_{|T}=t$. Since $\lambda\in F^{2}\Omega^{3} (K )$ by Corollary \ref{feb2211} we conclude that $d\omega(\check{t})_{|\pi^{-1}(b )}=0$ for all $b \in B$. Moreover, since $V$ is connected and $\omega(\check{t})_{|T}=0$ we conclude as at the end of the proof of Theorem \ref{feb1770} that
$[\omega(\check{t})_{|\pi^{-1}(b )}]=0\in H^{2}(\pi^{-1}(b);\R)$ for all $b \in V$.
Hence we can find a form $\beta\in \Omega^{1}(\pi^{-1}(V))$ such that $\omega(\check{t})_{|\pi^{-1}(b )}=d\beta_{|\pi^{-1}(b )}$. We now obtain a new section $\check{t}-\beta\in H^{{geom}}(\pi^{-1}(V))$ such that $\omega(\check{t}-\beta)=0$.
{Hence $\check{t}-\beta
\in H^{\flat}(\pi^{-1}(V))$ is the required section.}
\hB

%

The combination of the following Lemma with Lemma \ref{apr1101} and    Theorem \ref{feb1770}
clearly implies the second assertion of Theorem  \ref{feb2560}.

\begin{lem}\label{feb2450}
There exists a section $t\in H^{geom}(T)$
with $\omega(t)=0$ if and 
  only if the central extension \eqref{feb2430}
  is trivial. 
    \end{lem}
\proof
We use the fact that
there exists a section $t\in H^{geom}(T)$
with $\omega(t)=0$ if and 
  only if
  there exists a possibly different section $t\in H^{geom}(T)$ such that
$[\omega(t)]\in H^{2}(T;\R)$ is integral. 
In the following argument we show that such a section exists if and only if   the central extension \eqref{feb2430}
  is trivial. 
  We will further use that if $[\omega(t)]$ is integral for some section $t\in H^{geom}(T)$, then it is so for every section.

\bigskip

The idea is first to construct a section $\bar t\in H(T)$ and then discuss the obstruction to lift it to section $t\in H^{geom}(T)$ such that $[\omega(t)]\in H^{2}(T;\R)$ is integral. By the definition of a lifting gerbe in Section \ref{feb1302}, a section $\bar t\in H(T)$ is a $\widehat{LK}$-reduction of the $LK$-principal bundle $\tilde P_{|T}\to T$. In order to {construct} such a reduction we use the fact that a trivial $LK$ principal bundle has a canonical $\widehat{LK}$-reduction. Unfortunately $\tilde P_{|T}\to T$ is non-trivial, but its pull-back along the exponential map is so. In order to fix the section $\bar t$ we must fix a trivialization of this pull-back such that the associated canonical   $\widehat{LK}$-reduction   descents to $T$.


\bigskip


%

%

%



In detail this construction goes as follows.
We consider the composition 
$\phi: \taaa \stackrel{ \exp}{\to}  T  \to K  $ of the exponential map and the inclusion. The principal $K$-bundle  $(\phi\times \id_{S^{1}})^{*}P\to    \taaa\times S^{1}$ can be {trivialized}. In order to fix a {trivialization} we consider the
 section defined by $$s(u,[t]):=[\exp(tu),\exp(u) ,t]\ ,\quad  t\in \R\ ,\: u\in \taaa\ ,  $$
where $[t]\in S^{1}\cong \R/\Z$ denotes the class of $t\in \R$, and $[g,h,t]\in P$ denotes the point represented by $(h,g,t)\in K\times K\times \R$.
%
 The {trivialization} of ${(\phi\times \id_{S^{1}})^{*}}P$ induces a trivialization of the principal $LK$-bundle
 $\phi^{*}\tilde P\to   \taaa$. The corresponding section  is given by
 $$\tilde s(u)([t])=[\exp(tu),\exp(u),t]\ , \quad u\in \taaa\ , \:\:t\in \R\ .$$
 We let
 $$\check{P}\to \phi^{*}\tilde P$$ be the canonical $\widehat{LK}$-reduction associated to this trivialization.
 It has a section $\check s$ which lifts the section $\tilde s$.

 \bigskip
 We now show that the canonical $\widehat{LK}$-reduction descents to $T$.
 We have an action of $\Lambda$ on $\phi^{*}P$ given by
 $\lambda\cdot (u,[h,g,t])=(u+\lambda,[h,g,t])$ (note that $g=\exp(u)$). This action naturally induces an action of $\Lambda$ on  $\phi^{*}\tilde P$, and
 we have an identification  $\phi^{*}\tilde P/\Lambda\cong   \tilde P_{|T}$.
 
 \bigskip

The action of $\Lambda$ on $\phi^{*}P$  does not preserve the section $\tilde s$.
Indeed we have
$$\tilde s(u+\lambda)([t])=[\exp(t(u+\lambda)),\exp(u+\lambda),[t]]=(\tilde s(u)  \gamma_{\lambda})([t])\ .$$
In terms of the section $\tilde s$ the action of $\Lambda$ on $\phi^{*}\tilde P$  can be described 
in the form
\begin{equation}\label{feb2401}\lambda\cdot \tilde s(u)\ell=\tilde s(u+\lambda)\gamma_{\lambda}^{-1} \ell \ , \quad u\in \taaa\ ,\quad  \ell\in LK\ , \quad \lambda\in \Lambda\end{equation}
In a similar way we want to describe a lift of this action to $\check{P}$ using the section $\check{s}$.

\bigskip

   In the following we implicitly 
 use the natural embedding $\hat \Lambda\hookrightarrow \widehat{LK}$.
The  choice of a  set-theoretic split   $\tau_{0}:\Lambda\to \hat \Lambda$ determines 
 an antisymmetric bi-additive map (also called the commutator map)
\begin{equation}\label{apr1160}b:\Lambda\times \Lambda\to S^{1}\end{equation}
by \begin{equation}\label{feb2710}b(\lambda,\mu)=\tau_{0}(\lambda)^{-1}\tau_{0}(\mu)^{-1}\tau_{0}(\lambda)\tau_{0}(\mu)\ , \quad  \lambda,\mu\in \Lambda\ .\end{equation}
We can choose an {antisymmetric bilinear} form $B:\taaa\times \taaa\to \R$ such that 
\begin{equation}\label{feb2711}[B(\lambda,\mu)]_{\R/\Z}=b(\lambda,\mu)\ ,  \quad  \forall\lambda,\mu\in \Lambda\ .\end{equation} 
By the general theory of central extensions of free abelian groups of finite rank 
it is now possible to find another split
$\tau:\Lambda\to \hat \Lambda$ such that
\begin{equation}\label{feb2712}\tau(\lambda)\tau(\mu)=\tau(\lambda+\mu)[\frac{1}{2}B(\lambda,\mu)] \ ,  \quad  \forall\lambda,\mu\in \Lambda\ ,\end{equation}
i.e. $$\Lambda\times \Lambda\ni (\lambda,\mu)\mapsto [\frac{1}{2}B(\lambda,\mu)] \in S^{1}$$ is the cocycle of \eqref{apr1120} w.r.t the split $\tau$.
One can check that the following formula defines an action of $\Lambda$ on $\check{P}$ which covers the action on $\phi^{*}\tilde P$:
 \begin{equation}\label{feb2410}\lambda\cdot \check{s}(u)\hat \ell:= \check{s}(u+\lambda)\tau(\lambda)^{-1}  [\frac12 B(u,\lambda)]\hat\ell\ , \quad \lambda\in \Lambda , \quad  u\in \taaa\ , \quad\hat  \ell\in \widehat{LK}\ .\end{equation}
The quotient $\hat P:=\check{P}/\Lambda\to T$ is a $\widehat{LK}$-reduction of $\tilde P_{|T}\to T$
and therefore  a section $\bar t\in H(T)$.  

\bigskip

We now refine the section $\bar t$ to a section $t\in H^{geom}(T)$ by choosing a connection $\hat \alpha$ on $\hat P$ which lifts $\tilde \alpha$.
We start with the  connection one-form
$$\tilde s^{*}\tilde \alpha\in \Omega^{1}( \taaa)\otimes L\kaaa$$ for $\phi^{*}\tilde P$
pulled back by the section $\tilde s$. It follows from \eqref{feb2401}  that it satisfies
\begin{equation}\label{feb2411}\lambda^{*} \tilde s^{*}\tilde \alpha
=\Ad(\gamma_{\lambda}^{-1})\tilde s^{*}\tilde \alpha\ , \quad \lambda\in \Lambda\ .\end{equation}
We want to describe the connection $\hat \alpha$ on $\hat P$ in a similar way by providing the form
$\check{s}^{{*}}\hat{\alpha}\in \Omega^{1}(\taaa)\otimes \widehat{L\kaaa}$
which must satisfy $a(\check{s}^{{*}}\hat{\alpha})=\tilde s^{*}\tilde \alpha$ with $a$ as in \eqref{feb2402}. The analog of  the equivariance condition \eqref{feb2411} follows from \eqref{feb2410} and is given by
$$(\lambda^{*}\check{s}^{*}\hat{\alpha})(u)(X)=\Ad(\tau(\lambda)^{-1}) \check{s}^{*}\hat \alpha(u)(X)+ \frac12 B(X,\lambda) \ , \quad \lambda\in \Lambda\ ,\quad u,X\in \taaa\ .$$

These conditions are satisfied if we define $\hat \alpha$ such that
$$\check{s}^{*}\hat \alpha:=(-l_{\sigma}(\tilde s,\tilde s^{*}\tilde \alpha)+\alpha_{B}) \oplus_{\sigma} \tilde s^{*}\tilde \alpha\in \Omega^{1}( \taaa)\otimes \widehat{L\kaaa}\ ,$$
where $\alpha_{B}\in \Omega^{1}(\taaa)$ is given by $\alpha_{B}(u)(X)=\frac{1}{2}B(X,u)$.
 In this way we have defined a section $t=(\bar t,\hat \alpha)\in H^{geom}(T) $ which lifts $\bar t$.

\bigskip

We now calculate the curving  using using \eqref{feb1850}. We obtain $\omega(t)=B$, where we interpret $B$ as an invariant two-form on $T$.
We conclude that $[\omega(t)]$ is integral if and only if
$b=1$, i.e. if the central extension of $\Lambda$ is trivial. This finishes the proof of Lemma \ref{feb2450} and hence of Theorem \ref{feb2560}.
\hB
\begin{rem2}
In this remark we show that the condition of Theorem \ref{feb2560} might be non-trivial.
If $K$ is simply connected and simply laced, then it has been shown in \cite[4.8.1]{MR900587} that   the commutator map \eqref{apr1160} satisfies 
$$b=[\frac12 \langle\dots,\dots\rangle]_{\R/\Z}\ ,$$
where  $\langle\dots,\dots\rangle$ is the pairing \eqref{feb2440}.
If we take the basic central extension of a simply connected, simply laced group $K$ of rank $\ge 2$ like $SU(n)$, $n\ge 3$,
then the condition of Lemma \ref{feb2450}
 is not satisfied. In fact, since $K$ is simply connected, the integral lattice $\Lambda$ coincides with the coroot lattice, and  {there exist} two coroots
$H_{\alpha},H_{\beta}\in \Lambda$ with $\langle H_{\alpha},H_{\beta}\rangle=1$ and hence
$b(H_{\alpha},H_{\beta})=[\frac{1}{2}]\not=0$.

More {generally, the results} 
\cite[Theorem 3.2.1 and Theorem 3.5.1]{MR1724846} show that given an invariant symmetric bilinear form
$\langle\dots,\dots\rangle$ on $\kaaa$ and a map $b:\Lambda^{2}\Lambda\to S^{1}$   such that
$\langle H_{\alpha},H_{\beta}\rangle\in \Z$ for all pairs of fundamental coroots $H_{\alpha},H_{\beta}\in \taaa$ and   
$b(\lambda,H_{\alpha})=[\frac12 \langle \lambda,H_{\alpha}\rangle]$ for all $\lambda\in \Lambda$ and
coroots $H_{\alpha}$ there exists a central extension \eqref{feb2301}
yielding the map $b$ via \eqref{feb2710}.
This provides more examples.
\end{rem2}

\begin{rem2}
If the condition of Lemma \ref{feb2450} is not satisfied, then the canonical geometry on the lifting
gerbe $H\to K$ gives a reduction $H^{\flat}\to K$ whose restriction to the fibres of $K\to K/T$ is {nontrivial}. Of course, one could check using Theorem \ref{feb1901} that the set $\Red(H\to K)_{0}$ is not-empty. In view of Theorem \ref{feb2460} \eqref{pointtwo}, in order to get a reduction which is trivial on  the fibres of $K\to K/T$, one must modify the geometry by a two-form which is fibrewise closed but not exact.
Such a form cannot be closed since $H^{2}(G;\R)=0$. Hence this modification of the geometry necessarily changes the curvature which then can no longer be bi-invariant.
\end{rem2}

\begin{rem2}
The calculation above has an interpretation in terms of differential cohomology. In particular we can express the extension \eqref{feb2430} solely in terms of the class of $\beta([H\to K])\in H^{3}(K;\underline{\Z})$ and calculations in 
differential integral cohomology {as we will outline now.} 

\bigskip

 We consider a class $$x\in H^{3}(K;\underline{\Z})\ .$$
We have an  isomorphism
$$(\Lambda^{3}\kaaa^{*})^{K}\cong \Omega^{3}(K)^{K\times K}\cong H^{3}(K;\underline{\R^{\delta}})\ .$$
 There is a unique element
$c_x\in (\Lambda^{3}\kaaa^{*})^{K}$ such that the corresponding form $\lambda_x\in  \Omega^{3}(K)^{K\times K}$ represents the image of $x$ in $H^{3}(K;\underline{\R^{\delta}})$. There is a uniquely determined invariant pairing
$\langle\dots,\dots\rangle_x\in S^{2}\kaaa^{*}$ such that
\begin{equation}\label{feb2701}c_{x}(X,Y,Z)=\langle [X,Y],Z\rangle_{x}\ , \quad  X,Y,Z\in \kaaa \ .\end{equation}
 
The basic structure of differential cohomology  \cite{MR827262} provides
the exact {sequence}  
$$H^{2}(K;{\underline{\R^{\delta}}})\stackrel{a}{\to}\widehat{H\Z}^{3}(K)\stackrel{(I,R)}{\to}H^{3}(K;\underline{\Z})\times_{H^{3}(K;\underline{\R^{\delta}})}\Omega^{3}_{cl}(K)\to 0\ .$$
Since $H^{2}(K;\underline{\R^{\delta}})=0$ there exists a unique class $\hat x\in \widehat{H\Z}^{3}(K)$ such that $I(\hat x)=x$ and $R(\hat x)=\lambda_x$.

\bigskip

We consider a maximal torus $T\subseteq K$. By Lemma \ref{feb2211} we know that $\lambda_{x|T}=0$. We have an exact sequence
$$0\to H^{2}(T;{\bbG}^{\delta})\to \widehat{H\Z}^{3}(T)\stackrel{R}{\to} \Omega^{3}_{cl}(T)\ .$$
Therefore we can consider $\hat{x}_{|T}\in H^{2}(T;{\bbG}^{\delta})$.
The universal coeffficient theorem  gives an isomorphism
$$H^{2}(T;{\bbG}^{\delta})\cong \Hom(H_{2}(T;\underline{\Z});\Z)\cong \Hom(\Lambda^{2} \Lambda,S^{1})\ .$$
  We write
$$b_x:\Lambda^{2} \Lambda\to S^{1}$$ for the map determined by $\hat x_{|T}$.

\bigskip

We have a natural  isomorphism $$\Ext(\Lambda,S^{1})\cong \Hom(\Lambda^{2} \Lambda,S^{1})\ ,$$
where $\Ext(\Lambda,S^{1})$ is the group of central extensions of $\Lambda$ by $S^{1}$. In the following let us describe this identification. 
Given a central extension $\hat \Lambda$  we obtain the homomorphism $b:\Lambda^{2}\Lambda\to S^{1}$ by the formula \eqref{feb2710} using a set-theoretic split $\tau_0:\Lambda\to \hat \Lambda$.
On the other hand, given a map $b:\Lambda^{2}\Lambda\to S^{1}$ we choose a linear map $B:\Lambda^{2}\taaa\to \R$  satisfying \eqref{feb2711} and define a central extension $\hat \Lambda$ by defining the product on $S^{1}\times  
\Lambda$ by
$$(s,\gamma)\cdot (t,\mu)=(s+t+[\frac12 B(\gamma,\mu)],\gamma+\mu)\ , \quad s,t\in S^{1}\ , \quad \gamma,\mu\in \Lambda$$
(compare with \eqref{feb2712}).

\bigskip

The construction explained above thus associates to a class $x\in H^{3}(K;\underline{\Z})$ a
central extension class $b_{x}\in \Ext(\Lambda,S^{1})$.
 Recall that  central extension \eqref{feb2301} determines a lifting gerbe $H\to K$ as above and therefore a class $[H\to K]\in H^{2}(K;\bbG)$.

\begin{kor}
The extension class
$b_{\beta([H\to K])}\in \Ext(\Lambda,S^{1})$  is equal to the isomorphism class $b\in \Ext(\Lambda,S^{1})$ of the extension  \eqref{feb2430}.
\end{kor}
\proof
For $x:=\beta([H\to E])\in H^{3}(K,\underline{\Z})$ we have
$c_x=c$, $\lambda_x=\lambda$ and $\langle\dots,\dots\rangle_x=\langle\dots,\dots\rangle$,
where the left-hand sides are defined using the canonical geometry on $H\to K$ as above.
If we choose a geometry on the gerbe $H\to K$ with curvature $\lambda$, then the resulting geometric gerbe is classified by the lift $\hat x\in \widehat{H\Z}^{3}(K)$. Note that $H_{|T}\to T$ is flat and trivializable. The class
$\hat{x}_{|T}$ is the class of $[\omega(t)]\in H^{2}(T;\underline{\R^{{\delta}}})/H^{2}(T;\underline{\Z})\cong H^{2}(T;{\bbG}^{\delta})$ for every choice of section $t\in H^{geom}(T)$. In the proof of Lemma \ref{feb2450}   we have constructed a section
such that $\omega(t)\in \Omega^{2}_{cl}(T)^{T}$ is the invariant closed  form on $T$ determined by $B:\Lambda^{2}\taaa\to \R$ satisfying \eqref{feb2711}.
It follows that $b_{x}=b$. \hB

\end{rem2}

\section{Topological $T$-duality and Langlands duality}\label{feb2040}

The paper   \cite{2012arXiv1211.0763D} considers the torus bundle $K\to K/T$ for a compact Lie group $K$ and its maximal torus $T$ equipped with a twist $H\to K$, i.e. exactly the situation discussed in Section \ref{feb1920}.
It studies the question under which conditions it is $T$-dual to the corresponding bundle with twist for the Langlands dual group $K^{L}$.  

\bigskip

The  notion of   $T$-duality used  \cite[Definition 3.2]{2012arXiv1211.0763D} 
slightly differs from the integral definition of topological $T$-duality in terms of $T$-duality triple in \cite{MR2287642}. The concept in \cite[Definition 2.10]{2012arXiv1211.0763D}  follows \cite{MR2062361} and is based on differential forms. In the presence of torsion in the integral cohomology it may differ from 
our integral definition. 
 The definition of $T$-duality in \cite{2012arXiv1211.0763D}   is sufficient for the $T$-duality isomorphism in twisted cohomology, but in order to ensure a $T$-duality isomorphism in twisted $K$-theory one needs the integral version. 
 
 \bigskip
 
 Our goal in the present section is to provide an integral refinement 
   of the results in  \cite{2012arXiv1211.0763D}.
   In Theorem \ref{feb1904} we precisely describe the isomorphism classes of the gerbe $H\to K$ and its reduction $H^{\flat}\to K$
   which yields the Langlands dual group $K^{L}$ as $T$-dual.
   Apparently we get a more general result since we only have to exclude simple factors
of type B and C, while in \cite{2012arXiv1211.0763D} only factors of types ADE are {allowed}.

\bigskip

We consider a compact semisimple Lie group $K$ with maximal torus $T$. We exclude torus factors in order to simplify the presentation, but the theory extends smoothly to this more general case. 
The group $K$ is the total space of  a $T$-principal bundle over the flag variety $B:=K/T$.
First we collect some facts about the topology of the flag variety.

\begin{fact}
\begin{enumerate}
\item
The flag variety $B$ is simply connected.
\item $H^{*}(B;\Z)$ is torsion-free \cite{MR0071773}.
\item We have $H^{1}(B;\Z)=0$ and $H^{3}(B;\Z)=0$ \cite{2008arXiv0801.2444D}. 	
\end{enumerate}
\end{fact}
Let  $\tilde K\to K$ be the universal covering.  Then $\tilde T\to T $ denotes the corresponding covering of maximal tori. We add a symbol ${}^{\tilde{}}$ to all objects associated to $\tilde K$. 
There is a natural diffeomorphism
$B\cong \tilde K/\tilde T$.
  We have a natural identification
$\tilde \Lambda^{*}\cong H^{1}(\tilde T;\Z)$ where $\tilde \Lambda^{*}\subseteq \taaa^{*}$ denotes the lattice of integral weights.
We now consider the Leray-Serre spectral sequence $(\tilde E^{p,q}_r
,\tilde d_r)$ of the bundle $\tilde K\to B$. The differential of the second page  induces a map
$$b: \tilde \Lambda^{*}\cong H^{1}(\tilde T;\Z)\cong \tilde E_2^{0,1}\stackrel{\tilde d_2}{\to} \tilde E_2^{2,0}\cong H^{2}(B;\Z)\ .$$
It extends to a map of commutative algebras
$$b:S^{*}(\tilde \Lambda^{*})\to H^{even}(B;\Z)\ .$$

Let $S^{2}(\tilde \Lambda^{*})^{W}\subseteq S^{2}(\tilde \Lambda^{*})$ denote the subalgebra of Weyl invariant quadratic polynomials.

The following facts follow from the calculations in \cite{2008arXiv0801.2444D}.
\begin{fact}
\begin{enumerate}
\item  The map $b:\tilde \Lambda^{*}\to H^{2}(B;\Z)$ is an isomorphism.
\item  The map $b$ and the ring structure on $H^{even}(B;\Z)$ induce an isomorphism
$$ S^{2}(\tilde \Lambda^{*})/S^{2}(\tilde \Lambda^{*})^{W}\to H^{4}(B;\Z)\ .$$
\end{enumerate}
\end{fact}


%
%
%
%

We now consider the Leray-Serre spectral sequence $( E^{p,q}_r
,  d_r)$ of the bundle $  K\to B$. We have an identification
$\Lambda^{*}\cong H^{1}(T;\Z)\cong E_2^{0,1}$. 
By $F^{*}H^{*}(K;\Z)$ we denote the associated  decreasing filtration.
The edge sequence represents $H^{3}(K;\Z)$ as the cohomology of the complex
$$\Lambda^{2}( \Lambda^{*})\stackrel{d_2}{\to} \Lambda^{*}\otimes \tilde \Lambda^{*}\stackrel{d_2}{\to}  S^{2}(\tilde\Lambda^{*})/S^{2}(\tilde \Lambda^{*})^{W}\ ,$$
where the differentials are given by $d_2(\alpha \wedge \beta) =\beta\otimes  \iota^{*}(\alpha)-\alpha\otimes \iota^{*}(\beta)  $ and
$d_2(\alpha\otimes \beta)=\iota^{*}(\alpha)\beta$, and where we  use the map $\iota^{*}:\Lambda^{*}\to \tilde \Lambda^{*}$ which is dual to  $\iota:\tilde \Lambda\to \Lambda$.

\begin{kor}
 We have  $H^{3}(K;\Z)=F^{2} H^{3}(K;\Z)$. In particular,  all classes in $H^{3}(K;\Z)$
are $T$-dualizable in the sense of \cite{MR2287642}. \end{kor}

Equivalently, using the isomorphism
$H^{3}(K;\Z)\cong H^{3}(K;\underline{\Z} )$, for every gerbe $H\to K$ with band $\bbG$
the set $\Red(H\to K)_{0}$ is not empty by Corollary \ref{feb2601}.  For every choice of reduction in $\Red(H\to K)_0$ we get a  canonical $T$-dual $\hat H\to \hat K\to B$ as the right-hand side of the canonical $T$-duality diagram introduced in Definition \ref{feb1902}.

\bigskip

The classification of elements in $\Red(H\to E)_{0}$ is given in the second part of Theorem \ref{feb1901}. For the present purpose  we use the equivalent, but more explicit result  \cite[Theorem 2.24]{MR2287642}. As above we identify $\Lambda^{*}\otimes \tilde \Lambda^{*}\cong \Hom(\Lambda,\tilde \Lambda^{*})\cong E_{2}^{2,1}$. Let $h\in H^{3}(K;\Z)$ be the class of a gerbe $H\to K$ with band $\bbG$. In the following we specialize \cite[Theorem 2.24]{MR2287642} using in particular the fact that  $H^{3}(B;\Z)=0$.

\begin{theorem}[\cite{MR2287642}]\label{feb2001}
The extension of the pair $(K\to B,h\in H^{3}(K;\Z))$ to a $T$-duality triple, 
 or equivalently, an element in $\Red(H\to K)_{0}$,
   is determined by the choice of a representative
$u\in \Hom(\Lambda,\tilde \Lambda^{*})$ of the class $h$.
 In particular,
the Chern classes of the $T$-dual $T$-bundle $\hat K\to B$ are the generators of the image $u(\Lambda)\subseteq \tilde \Lambda^{*}$ of $u$.
 \end{theorem}

We refer to \cite[Section 2]{2012arXiv1211.0763D} for a short introduction to Langlands duality.
By $\kaaa^{L}$ we denote the Langlands dual of the Lie algebra $\kaaa$. We now assume that no simple factor of $\kaaa$  is of  type BC. Then there exists
an isomorphism of Lie algebras $\phi:\kaaa\to \kaaa^{L}$. 
Its restriction to $\taaa$  preserves integrality. Consequently, $\phi$ induces an isomorphism
$\phi_{|\tilde \Lambda^{L,*}}^{*}:\tilde \Lambda^{L,*}\to \tilde \Lambda^{*}$. Langlands duality yields an inclusion
$\lambda:\Lambda \cong \Lambda^{L,*}\to \tilde \Lambda^{L,*}$. We define
 \begin{equation}\label{nov3100}u:=\phi\circ \lambda:\Lambda\to \tilde \Lambda^{*}\ . \end{equation}
\begin{theorem}\label{feb1904}
The map defined in \eqref{nov3100} is a cycle. The total space of the associated $T$-dual
is the Langlands dual group $K^{L}$.
\end{theorem}
\proof
Note that $d_2(u)\in S^{2}(\tilde \Lambda^{*})/S^{2}(\tilde \Lambda^{*})^{W}$
is represented by the bilinear form $(X,Y)\mapsto u(X)(Y)$, $X,Y\in \tilde \Lambda$.
Since $u$ is Weyl invariant, this bilinear form is Weyl invariant, too. Hence $u$ is a cycle.

\bigskip

The map $\phi$ induces a diffeomorphism of flag varieties $B(\phi): B=B(\kaaa)\to B(\kaaa^{L})$ which we use to view $K^{L}$ as a $T^{L}$-principal bundle over $B$.
Up to diffeomorphism the total space of the torus bundle $K^{L}\to  B$ is determined by the image of $d_2^{L}:E_2^{L,0,1}\to E_2^{L,2,0}$.  As before we identify this map with     the map $\Lambda^{L,*}\hookrightarrow \tilde \Lambda^{L,*}$, where $(E^{L,p,q}_r,d_r^{L})$ is the Leray-Serre spectral sequence of $K^{L}\to B$. Pulling back by $\phi$ we get the inclusion of lattices  $\phi_{|\taaa}^{*}( \Lambda^{L,*})\hookrightarrow \tilde \Lambda^{ *}$. Now note that $\Lambda^{L,*}=\lambda(\Lambda)$.  \hB

The observation that one can get the Langlands dual group $K^{L}$  of $K$ by $T$-duality is the basis for more interesting question. For example, how can one describe the group structure on $K^{L}$ in terms of data on $K$ if one realizes $K^{L}$ as the moduli stack of fibrewise trivialization of a reduced gerbe
$H^{\flat}\to K$?

\bigskip

\section{A technical result}

Let $\cG:G^{1}\Rightarrow G^{0}$ be a   Lie groupoid which acts properly on a vector bundle $V\to M$.  Then we get a morphism of action groupoids
$$(G^{1}\times_{G^{0}}V\Rightarrow V)\to (G^{1}\times_{G^{0}}M\Rightarrow M)\ .$$
{Taking} nerves we obtain the simplicial vector   bundle
$$V^{\bullet}\to M^{\bullet}\ .$$
 We denote a point in $M^{p}$ by $(\gamma_{1},\dots,\gamma_{p-1},m)$. 
The face maps $\partial_{i}$ of $M^{\bullet}$ are given by
$$(\gamma_{1},\dots,\gamma_{p},m)\mapsto \left\{ \begin{array}{cc} (\gamma_{2},\dots,\gamma_{p},m)&i=0\\
(\gamma_{1},\dots, \gamma_{i}\gamma_{i+1},\dots,\gamma_{p},m)&i=1,\dots,p-1\\
(\gamma_{1},\dots,\gamma_{p}m)&i=p
 \end{array}\right. \ ,$$
and similarly for $V^{\bullet}$.
From now we will implicitly use the canonical identification of  fibres $V^p_{\gamma_1,\dots,\gamma_p,m} \cong V_m$ for 
$(\gamma_1,\dots,\gamma_p,m) \in M^p$. Thus for a section  $f:M^{p}\to V^{p}$ we have
$$f(\gamma_{1},\dots,\gamma_{p},m)\in V_{m}\ .$$ 
This identification will be important for the correct interpretation of formula \eqref{glmai5}.
We set $d_{i}f:=\partial_{i}^{*} f$ where the pullback of sections involves the groupoid action on $V$. In detail we have
$$(d_{i}f)(\gamma_{1},\dots,\gamma_{p+1},m)= \left\{ \begin{array}{cc} f(\gamma_{2},\dots,\gamma_{p+1},m)&i=0\\
f(\gamma_{1},\dots, \gamma_{i}\gamma_{i+1},\dots,\gamma_{p+1},m)&i=1,\dots,p\\
\gamma_{p+1}^{-1} f(\gamma_{1},\dots,\gamma_{p},\gamma_{p+1}m)&i=p+1
\end{array}\right.\ .$$
We consider the complex
$$(C^{{\infty}}(M^{\bullet},V^{\bullet}),d)\ , \quad  d:=\sum_{i=0}^{p+1} (-1)^{i}d_{i}:C^{{\infty}}(M^{p},V^{p})\to C^{{\infty}}(M^{p+1},V^{p+1})\ .$$
\begin{lem}\label{apr2301}
We have  $H^{p}(C^{{\infty}}(M^{\bullet},V^{\bullet}),d)=0$ for $p\ge 1$.
\end{lem}
\proof  
Since the  groupoid $G^{1}\Rightarrow^{t}_{s} G^{0}$ is locally compact we can choose a smooth Haar system $(\lambda^{g})_{g\in G^{0}}$. Here $\lambda^{g}$ is a measure on $G^{g}:=t^{-1}(g)$ such that $\gamma^{-1}_{*}\lambda^{t(\gamma)}=\lambda^{s(\gamma)}$ for all $\gamma\in  G^{1}$.
 
 Since  the action groupoid $G^{1}\times_{G^{0}}M\Rightarrow^{\mu}_{\pr} M$ is proper  by \cite[Prop. 6.11]{MR1671260}
there exists a function
$\chi\in C^{\infty}(M)$ such that $t: \supp( \mu^{*} \chi)\to G^{0}$ is proper and  $$\int_{G^{r(m)}} \chi(\gamma^{-1} m) d\lambda^{r(m)}(\gamma)=1$$
for all $m\in M$, where $r:M\to G^{0}$ is the structure map.

\bigskip

Let $p\ge 0$ and
 $f\in C^{{\infty}}(M^{p+1},V^{p+1})$ be a cycle. Then we define
$h\in C^{{\infty}}(M^{p},V^{p})$ by
\begin{equation}\label{glmai5}
h(\gamma_{1},\dots,\gamma_{p-1},m):=\int_{G^{s(\gamma_{p-1})} } \chi(  \gamma^{-1}  m) \gamma f(\gamma_{1},\dots,\gamma_{p-1},\gamma, \gamma^{-1}m) d\lambda^{s(\gamma_{p-1})}(\gamma)\ .
\end{equation}
Since $f$ is a cycle we have 
\begin{eqnarray*}0&=&\gamma f(\gamma_{2},\dots,\gamma_{p},\gamma,\gamma^{-1}m)\\
&&+\sum_{i=1}^{p-1}(-1)^{i} \gamma f(\gamma_{1},\dots,\gamma_{i}\gamma_{i+1},\dots,\gamma_{p},\gamma,\gamma^{-1}m)\\ &&+  (-1)^{p}\gamma f(\gamma_{1},\dots, \gamma_{p-1},\gamma_{p}\gamma,\gamma^{-1}m)+ (-1)^{p+1}  f(\gamma_{1},\dots,\gamma_{p-1}, m)\ .
\end{eqnarray*}
Employing this formula in the equality marked by $!!$ and the invariance of the Haar system in the equality marked by $!$ we calculate
\begin{eqnarray*}
\lefteqn{(dh)(\gamma_{1},\dots,\gamma_{p},m)}&&\\
&=&h(\gamma_{2},\dots,\gamma_{p},m)\\
&&+
\sum_{i=1}^{p-1} (-1)^{i}h(\gamma_{1},\dots,\gamma_{i}\gamma_{i+1},\dots,\gamma_{p},m)\\
&&+(-1)^{p }\gamma_{p}^{-1} h(\gamma_{1},\dots, \dots,\gamma_{p-1},\gamma_{p}m)\\
&=&\int_{G^{s(\gamma_{p})}} \chi(\gamma^{-1}m) \gamma f(\gamma_{2},\dots,\gamma_{i}\gamma_{i+1},\dots,\gamma_{p},\gamma,\gamma^{-1}m)   d\lambda^{s(\gamma_{p})}(\gamma)
\\&&+
\sum_{i=1}^{p-1} (-1)^{i} \int_{G^{s(\gamma_{p})}} \chi(\gamma^{-1}m) \gamma f(\gamma_{1},\dots,\gamma_{i}\gamma_{i+1},\dots,\gamma_{p},\gamma,\gamma^{-1}m)   d\lambda^{s(\gamma_{p})}(\gamma)\\
&&+(-1)^{p}\int_{G^{s(\gamma_{p-1})}}\chi(\gamma^{-1}\gamma_{p}m)  \gamma_{p}^{-1}  \gamma f(\gamma_{1},\dots, \dots,\gamma_{p-1},\gamma,\gamma^{-1}\gamma_{p}m)\lambda^{s(\gamma_{p-1})}(\gamma)\\&\stackrel{!}{=}&\int_{G^{s(\gamma_{p})}} \chi(\gamma^{-1}m) \gamma f(\gamma_{2},\dots,\gamma_{i}\gamma_{i+1},\dots,\gamma_{p},\gamma,\gamma^{-1}m)   d\lambda^{s(\gamma_{p})}(\gamma)
\\&&+\sum_{i=0}^{p-1} (-1)^{i} \int_{G^{s(\gamma_{p})}} \chi(\gamma^{-1}m) \gamma f(\gamma_{1},\dots,\gamma_{i}\gamma_{i+1},\dots,\gamma_{p},\gamma,\gamma^{-1}m)   d\lambda^{s(\gamma_{p})}(\gamma)\\
&&+(-1)^{p}\int_{G^{s(\gamma_{p})}}\chi(\gamma^{-1}m)     \gamma f(\gamma_{1},\dots, \dots,\gamma_{p-1},\gamma_{p}\gamma,\gamma^{-1}m)\lambda^{s(\gamma_{p})}(\gamma)\\
&\stackrel{!!}{=}& (-1)^{p} \int_{G^{s(\gamma_{p})}}\chi(\gamma^{-1}m)   f(\gamma_{1},\dots,\gamma_{p},m) \lambda^{s(\gamma_{p})}(\gamma)\\
 &=&(-1)^{p} f(\gamma_{1},\dots,\gamma_{p},m) \ ,
\end{eqnarray*}
i.e. $$d h=(-1)^{p}f\ .$$\hB  
\section{Conventions}

Since some of the calculations are sensitive for the choice of normalizations in the Cartan calculus in the following we list the conventions used in the present paper.
\begin{enumerate}
\item $(d\omega)(X,Y)=Y\omega(X)-Y\omega(X)$ for a one form $\omega$ and   vector fields $X,Y$ with $[X,Y]=0$
\item $(d\omega)(X,Y,Z)=X\omega(Y,Z)-Y\omega(X,Z)+Z\omega(X,Y)$ for a three form $\omega$ and vector fields  $X,Y,Z$ with $[X,Y]=[Y,Z]=[X,Z]=0$
\item $(\alpha\wedge \beta)(X,Y)=\alpha(X)\beta(Y)-\alpha(Y)\beta(X)$  for two one forms $\alpha,\beta$ and vector fields $X,Y$.
\item The curvature of a connection $\alpha$ on a principal bundle  is given by $\Omega=d\alpha+[\alpha\wedge \alpha]$, i.e.
$$\Omega(X,Y)=X\alpha(Y)-Y\alpha(X)+2[\alpha(X),\alpha(Y)]$$ 
for commuting vector fields $X,Y$.
\item We identify $S^{1}\cong \R/\Z$. This induces the identification of its Lie algebra with $\R$.
\end{enumerate}

\bibliographystyle{alpha}
\bibliography{tred}

\end{document}